\documentclass[12pt]{article}
\usepackage{amsfonts,amssymb,amsmath,amsthm,eucal}
\numberwithin{equation}{section}

\newcommand{\la}{\lambda}
\newcommand{\al}{\alpha}
\newcommand{\Ls}{\Lambda^*}
\newcommand{\Z}{\mathbb{Z}}
\newcommand{\fz}{\mathfrak{z}}
\newcommand{\R}{\mathbb{R}}
\newcommand{\C}{\mathbb{C}}
\newcommand{\cD}{\mathcal{D}}
\newcommand{\Covc}{\CC} 
\newcommand{\e}[1]{\mathbf{e}^{#1}} 
\newcommand{\Th}{\vartheta} 
\newcommand{\tA}{\widetilde{A}} 
\newcommand{\tF}{\widetilde{F}} 
\newcommand{\cS}{\mathcal{S}} 
\newcommand{\cT}{\mathcal{T}}

\newcommand{\Aq}{A^{[\textup{?}]}} 
\newcommand{\Al}{A_{\textup{lead}}} 
\newcommand{\hn}{\widehat{n}} 
\newcommand{\fm}{\mathfrak{M}} 
\newcommand{\fP}{\mathfrak{P}} 
\newcommand{\tip}{\widetilde{p}}

\newcommand{\vo}{\mathbf{c}}
\newcommand{\tp}{\,\top\,}
\newcommand{\vk}{\varkappa}
\newcommand{\CC}{\mathcal{C}}

\DeclareMathOperator{\Cov}{Cov}
\DeclareMathOperator{\Covar}{Covar}
\DeclareMathOperator{\Aut}{Aut}
\DeclareMathOperator{\Hom}{Hom}
\DeclareMathOperator{\wt}{wt}

\DeclareMathOperator{\sgn}{sgn}

\newcommand{\cA}{\mathbf{A}}

\newcommand{\hb}{h} 

\newcommand{\lan}{\left\langle}
\newcommand{\llan}{\lan\!\lan\,}
\newcommand{\ran}{\right\rangle}
\newcommand{\rran}{\,\ran\!\ran}
\newcommand{\lran}{\lan\,\cdot\,\ran}
\newcommand{\lal}{\lan}
\newcommand{\ral}{\ran_{\hb}}

\newtheorem{theorem}{Theorem}[section]
\newtheorem{lemma}[theorem]{Lemma}
\newtheorem{proposition}[theorem]{Proposition}

\theoremstyle{definition}
\newtheorem{remark}[theorem]{Remark}
\newtheorem{definition}[theorem]{Definition}
\newtheorem{example}[theorem]{Example}

\newcommand{\bb}{\mathbb}

\newcommand{\cx}{{\bb C}}

\newcommand{\integers}{{\bb Z}}

\newcommand{\ratls}{{\bb Q}}
\newcommand{\reals}{{\bb R}}

\newcommand{\isom}{\cong}
\newcommand{\area}{\operatorname{Area}}
\newcommand{\zed}{\integers}
\newcommand{\vol}{\operatorname{vol}}

\newcommand{\cH}{{\cal H}}

\title{Asymptotics of numbers of branched coverings of a torus
and volumes 
  of moduli spaces of holomorphic differentials}
\author{Alex Eskin\thanks{
Department of Mathematics, University of Chicago,
5734 University Ave., Chicago, IL 60637.
E-mail: eskin@math.uchicago.edu}
{}  and Andrei Okounkov\thanks{
 Department of Mathematics, University of California at
Berkeley, Evans Hall \#3840, 
Berkeley, CA 94720-3840. E-mail: okounkov@math.berkeley.edu}}
\date{}

\begin{document}
\maketitle

\tableofcontents

\section{Introduction}

\subsection{Moduli spaces of holomorphic differentials}

Let $\Sigma$ be  a compact Riemann
surface of genus $g>1$ and $\omega$ is a holomorphic 1-form on $\Sigma$, 
i.e.\ a tensor of the form
$\omega(z)dz$ in local coordinates with $\omega$ holomorphic.  Away from
the zeros of $\omega$, we can
choose a  coordinate $\zeta$ so that
$\phi=d\zeta$. This determines a Euclidean metric $|d\zeta^2|$ in
that chart and the coordinate changes
between such charts 
are of the form $\zeta\to \zeta +c$.  
Consequently, holomorphic differentials are sometimes referred
to as translation surfaces or flat structures with parallel vector
fields.  

Near a zero of order $k\geq 1$ of $\omega$, we can choose
a local coordinate $\zeta$  so that $\omega$ 
is given by $\zeta^k \, d\zeta$. The corresponding
metric is then $|\zeta^{2k}||
d\zeta^2|$.  The total angle around the zero is $(2k+2)\pi$, so we say
that $\Sigma$ has a cone singularity with total angle $(2k+2)\pi$. 

\begin{definition}
Suppose that $g>1$ and let $\mu$ be a partition of $2g-2$ into
$\ell=\ell(\mu)$ parts. We denote by
$\cH(\mu)$ the moduli space of $(\ell+2)$-tuples
$$
(\Sigma,\omega,p_1,\dots,p_\ell)\,,
$$ 
where $\Sigma$ is a  Riemann surface of genus $g$, and $\omega$ is
an holomorphic differential on $\Sigma$, and
$$
(\omega)=\sum_i \mu_i \, [p_i] \,,
$$
where $(\omega)$ is the divisor of $\omega$, that is, the set
of zeros of $\omega$ counting multiplicity. 
\end{definition}

For example if $\mu = (3,1)$, 
we require that $\omega$ has one triple zero $p_1$ and one simple zero $p_2$.
Similarly, one can consider moduli spaces of 
pairs $(\Sigma,\omega)$ without ordering of the zeros of $\omega$. This
presents no additional difficulties. 

One important feature is that these spaces, and the similar spaces of 
quadratic differentials, admit 
an ergodic $SL(2,\reals)$ action.  The dynamics of this action is related to 
billiards in rational polygons and to interval exchange transformations. 
This circle of ideas has been studied extensively by various authors, 
e.g. \cite{Masur:interval, Masur:Smilley:hausdorff,
Veech:geodesic, Veech:moduli, Veech:Siegel,
Kontsevich:Zorich:components, Kra:Nielsen, Veech:surfaces,
Earle:Gardiner:Teichmuller, Gutkin:Judge}.

\subsection{Local coordinates and invariant measure on $\cH(\mu)$}

Consider the relative homology group 
$$
H_1(\Sigma,\{p_i\},\zed)\isom \Z^n\,, \quad 
n = 2g+\ell(\mu)-1\,, 
$$
where $\ell(\mu)$ is the number of parts in the partition $\mu$. 
Choose a basis
$$
\{ \gamma_1, \dots, \gamma_n \} \subset H_1(\Sigma,\{p_i\},\zed)
$$ 
so that $\gamma_i$, $i=1,\dots,2g$, form a standard symplectic basis
of $H_1(\Sigma,\zed)$ and 
$$
\partial\gamma_{2g+i}=[p_{i+1}]-[p_1]\,,\quad i=1,\dots,\ell(\mu)-1\,.
$$
The group $Sp(2g,\Z)\ltimes \Z^{2g(\ell-1)}$ acts transitively on such
bases by changing the basis in $H_1(\Sigma,\zed)$ and translating the cycles
$\gamma_{2g+i}$ by elements of $H_1(\Sigma,\zed)$. Consider the period map 
$$
\Phi: \cH(\mu) \to \cx^n \isom
(\reals^2)^n
$$
defined by
\begin{displaymath}
\Phi(\Sigma,\omega) = \left( \int_{\gamma_1} \omega, \dots, \int_{\gamma_n}
  \omega \right)
\end{displaymath}
It is known \cite{Kontsevich:lyapunov}  
 that $\Phi$ is a local
coordinate system on $\cH(\mu)$. In particular,
\begin{equation}\label{dimH}
\dim_\cx \cH(\mu) = 2g + \ell(\mu)-1\,.
\end{equation}

Let us pull back the Lebesgue measure from $\C^n$ to $\cH(\mu)$ using $\Phi$. 
This is well defined
since it is clearly independent of the choice of basis $\{\gamma_i\}$. 
This measure is infinite, essentially because $\omega$ can be
multiplied by any complex number. To correct this, we introduce
the following

\begin{definition} Denote by $\cH_1(\mu)$ the subset of $\cH(\mu)$
defined by the equation $\area_\omega(\Sigma) = 1$, where
$$
\area_\omega(\Sigma) = \frac{\sqrt{-1}}{2} 
\int_\Sigma \omega \wedge \overline{\omega}
$$
 is the area of $\Sigma$ with respect to the metric
defined by $\omega$. 
\end{definition}

In terms of the periods $\phi=\Phi(\Sigma,\omega)$ we have 
\begin{equation}\label{defQ}
\area_\omega(\Sigma) = 
\frac12 \sum_{i=1}^g \left(\phi_i \bar\phi_{g+i} -
\bar\phi_{i} \phi_{g+i} \right)
\end{equation}
Denote by $Q$ the quadratic form on $\R^{2\dim\cH(\mu)}$ defined by
\eqref{defQ}. It follows that the image of $\cH_1(\mu)$ under 
$\Phi$ is contained in the
hyperboloid $Q(v)=1$ and $\cH_1(\mu)$
can be identified with a certain open subset of $Q(v)=1$. 
We now define a measure $\nu$ on $\cH_1(\mu)$ as follows. 

\begin{definition}
Let a set $E \subset \cH_1(\mu)$ lie in the domain
of a coordinate chart $\Phi$ and let $C\Phi(E)\subset \C^n$ be the 
cone over $\Phi(E)$ with vertex at the origin $0\in\C^n$. 
By definition, we set
$$
\nu(E) = \vol(C\Phi(E)) \,,
$$
where the volume on the right is with respect to the Lebesgue
measure on $\C^n$. 
\end{definition}

This measure is invariant under the $SL(2,\reals)$ action on
$\cH_1(\mu)$. It is a theorem of Masur \cite{Masur:interval}  and Veech 
 \cite{Veech:Gauss} that 
$$
\nu(\cH_1(\mu)) < \infty\,.
$$ 
The main goal of this paper is the computation of these numbers. They
arise in particular in problems associated with billiards in rational
polygons \cite{Veech:Siegel, Eskin:Masur:ae, Eskin:Masur:Zorich:SV},
 and also in connection with interval
exchanges and the Lyapunov exponents of the Teichmuller geodesic flow
\cite{Zorich:interval, Kontsevich:lyapunov}.

\subsection{Volumes and branched coverings}

Our approach to the computation of the numbers $\nu(\cH_1(\mu))$
is based on the interpretation of $\nu(\cH_1(\mu))$ as the
asymptotics in a certain enumeration problem, namely,
the enumeration of connected branched coverings of
a torus as their degree goes to $\infty$ and the
ramification type is fixed. This interpretation was discovered by 
Kontsevich and Zorich and, independently, by Masur and the first author.  

\begin{definition}
Given a partition $\mu$, denote by $\CC_d(\mu)$ the weighted
the number of connected ramified coverings of
the standard torus 
\begin{equation}\label{siMT}
\sigma: \Sigma \to T
\end{equation}
of degree $d$, which are ramified over $\ell(\mu)$ fixed points
of $T$, and such that the nontrivial part of the monodromy around the
$i$th point is a cycle of length $\mu_i$. The weight of a covering
\eqref{siMT} is $\left|\Aut(\sigma)\right|^{-1}$,
where $\Aut(\sigma)$ is the commutant of the monodromy
subgroup of $\sigma$ inside the symmetric group $S(d)$. 
\end{definition}

\begin{remark} Typically, the group $\Aut(\sigma)$ is trivial
and, in particular, these weights make no impact on the asymptotics
of $\CC_d(\mu)$ as $d\to\infty$, see Section \ref{auto}. The purpose
of introducing the weights is to make certain exact formulas
look better, such as, for example, to make the generating series \eqref{Cq} 
a quasimodular form. 
\end{remark}
 
\begin{proposition}
\label{prop:count:tori:and:volumes} For any partition $\mu$, we have 
\begin{equation}\label{volcov}
\nu(\cH_1(\mu)) = \lim_{D\to\infty} D^{-\dim_\C \cH(\mu)} \, \sum_{d=1}^D \CC_d(\mu+\vec1) \,,
\end{equation}
where $\mu+\vec 1 =(\mu_1+1,\dots,\mu_{\ell(\mu)}+1)$.  
\end{proposition}

Recall that the dimension of $\cH(\mu)$ is given by \eqref{dimH}. 
The proof of Proposition \ref{prop:count:tori:and:volumes} is
elementary and is supplied in Section \ref{pofp} below.
The basic idea behind Proposition \ref{prop:count:tori:and:volumes} is that
to any covering \eqref{siMT} we can associate the
point 
$$
\left(\Sigma, \sigma^*(dz)\right)\in \cH(\mu)\,,
$$
where $dz$ is the standard holomorphic
differential on $T$, and counting such points in $\cH$ is like
counting points of $\Z^{2n}$ inside subsets of $\R^{2n}$,
where $n=\dim \cH(\mu)$.  

Using Proposition~\ref{prop:count:tori:and:volumes},  A.~Zorich computed
the numbers $\nu(\cH_1(\mu))$ for small $\mu$.

\subsection{Enumeration of coverings}

It will be convenient to introduce the following numbers
\begin{equation}\label{defvo}
\vo(\mu) = (|\mu|+1)\,  \lim_{D\to\infty}  D^{-|\mu|-1} \, 
\sum_{d=1}^D \CC_d(\mu)\,,
\end{equation}
where $|\mu|=\sum \mu_i$. The existence of this limit follows from
Proposition \ref{prop:count:tori:and:volumes} which states that 
$$
\vol(\cH_1(\mu)) = \frac{\vo(\mu+\vec 1)}{\dim \cH(\mu)} \,.
$$ 
Heuristically, one should think about \eqref{defvo} as
saying that
$$
\CC_d(\mu) \approx \vo(\mu) \, d^{|\mu|} 
$$
for a typical large number $d$. In this paper, we obtain a general 
formula for these numbers, and hence for the volumes $\nu(\cH_1)$, 
by developing a systematic approach to the
asymptotics of $\CC_d(\mu)$ as $d\to\infty$.  

Our starting point is an exact result of S.~Bloch and the
second author \cite{BO}, who considered certain 
generating functions, called the $n$-point functions,
which encode the numbers $\CC_d(\mu)$. These
$n$-point functions were evaluated in \cite{BO}
in a closed form as a determinant of $\vartheta$-functions 
and their derivatives, see also the paper \cite{O} for
a simplified approach. This result is
reproduced in Theorem \ref{t1}  below.
A qualitative conclusion from it
is that the following generating function 
\begin{equation}\label{Cq}
\CC(\mu) = \sum_{d=0}^\infty q^d \, \CC_d(\mu)
\end{equation}
is a \emph{quasimodular form} in the variable $q$ for the full modular
groups, that is, a polynomial in the 
Eisenstein series $G_k(q)$, $k=2,4,6$.  

The asymptotics in \eqref{defvo} corresponds to the
$q\to 1$ asymptotics of \eqref{Cq}. In principle,
using the formula for the $n$-point function, one can express
for any given $\mu$ the generating function \eqref{Cq} in Eisenstein
series. The quasimodularity of $\CC(\mu)$ means that
it transforms in a certain way under the transformation
$$
q=e^{2\pi i \tau} \mapsto e^{-2\pi i/\tau}\,,
$$
which takes $q=1$ to $q=0$, thus giving the $q\to 1$
asymptotics of $\CC(\mu)$. This quasimodularity is
a manifestation of a certain ``mirror symmetry'' between
coverings of very large degree ($q\to 1$) and small
degree ($q\to 0$). 

In practice, however, the computation of \eqref{Cq}
becomes very difficult even for relatively small $\mu$.

We therefore pursue a different approach and first
investigate the $q\to 1$ asymptotics of the $n$-point 
function. Here we find a great simplification, see
Theorem \ref{t2}, essentially because the
$\vartheta$-functions become trigonometric
functions. We then extract from this asymptotics the
information about the asymptotics \eqref{defvo}.

This extraction is still rather nontrivial because, 
inside of the $n$-point function,
the numbers $\CC(\mu)$ are  wrapped up in several layers
of enciphering, such as going from connected to disconnected coverings.
For example, the terms corresponding to connected coverings
appear only deep in the asymptotic expansion of the $n$-point
function which requires us to keep track of many orders of
asymptotics.  

\subsection{Summary of results} 

The answer we obtain for the constants $\vo(\mu)$ can be conveniently
stated  in terms of a certain multilinear form 
\begin{equation}\label{lala}
\lal \cdot | \dots | \cdot \ral :
\Ls \times \cdots \times \Ls \to \C[\hb^{-1}]\,,
\end{equation}
where $\Lambda^*$ an algebra closely related to the
algebra of symmetric functions.  The form \eqref{lala} is
such that 
\begin{equation}\label{vmu}
\lal f_{\mu_1} | \dots | f_{\mu_k} \ral = \vo(\mu) \, 
\frac{|\mu|!}{\hb^{|\mu|+1}}+\dots\,,
\end{equation}
where $f_k$ are certain generators of $\Lambda^*$, and dots
stand for terms of lower degree in $\hb^{-1}$. The evaluation of 
\eqref{vmu} goes in 3 steps. 

First, one expresses the
generators $f_k$ as polynomials in power-sum generators $p_k$
of $\Lambda$. A  formula for this
expansion is obtained in Theorem \ref{t4}. Then, using
an analog of the Wick formula for \eqref{lala}  
derived in Theorem \ref{t5}, one reduces \eqref{vmu} to
computations of the constants $\llan \mu \rran$ 
defined by 
$$
\lal p_{\mu_1} | \dots | p_{\mu_k} \ral = 
\frac{\llan \mu \rran}{\hb^{|\mu|+1}}+\dots\,,
$$
which we call elementary cumulants. 

These numbers $\llan \mu \rran$ 
are finally computed in Theorem \ref{tem}
in terms of values of the $\zeta$-function at even
positive integers, that is, in Bernoulli numbers.
In particular, we have 
$$
\pi^{-2g} \, \nu(\cH_1(\mu)) \in \mathbb{Q}
$$
for any $\mu$. This rationality was also conjectured by
Kontsevich and Zorich.

We were unable to simplify this answer further in
the general case, but in the special case $\mu=(2,\dots,2)$ 
which corresponds to differentials with simple
zeros (that is, to generic ones), 
an attractive answer is available. It is given in 
Theorem \ref{t7}.

\subsection{Example of a volume computation}

Suppose we want to compute $\nu(\cH(3,1))$ or, equivalently, 
$\vo(4,2)$. From Theorem  \ref{t4}
we get
$$
f_2 = \frac12\, p_2\,, \quad f_4 = \frac14\, p_4 - p_2 p_1 + \dots \,, 
$$
where dots stand for lower weight term which make
no contribution to the answer.

In general, there exist a very important weight filtration
on $\Lambda^*$ which we discuss in Section \ref{s4}. 
It has the property that \eqref{lala}  takes it to
the filtration of $\C[\hb^{-1}]$ by degree, which
allows us to identify many negligible terms.

By the Wick formula, see Theorem \ref{t5}, we have
\begin{multline*} 
\lal f_4 \big | f_2 \ral = 
\frac18 \lal p_4 \big | p_2 \ral - 
\frac12 \lal p_2 p_1 \big | p_2 \ral +\dots =\\
\frac{\hb^{-7}}{8} \llan 4,2 \rran - 
 \frac{\hb^{-7}}{2} \llan 2  \rran \, \llan 2,1 \rran -
\frac{\hb^{-7}}{2} \llan 1  \rran \, \llan 2,2 \rran + \dots \,,
\end{multline*}
where dots stand for lower terms. From Theorem \ref{tem},
see also Example \ref{covar}, we conclude that
$$
 \llan 1  \rran =  \zeta(2) = \frac{\pi^2}{6} \,, \quad
 \llan 2  \rran  = 0 \,,
$$
and similarly
$$
\llan 4,2 \rran =\frac{416}{315} \, \pi^6\,, \quad
\llan 2,2 \rran =\frac{16}{45} \, \pi^4\,.
$$
Hence
$$
\lal f_4 \big | f_2 \ral = \frac{128}{945} \, \pi^6\, \hb^{-7} + \dots \,,
$$
which means that
$$
\vo(4,2) = \frac{8}{42525}\, \pi^6 \,, \quad
\nu(\cH(3,1))=\frac{8}{297675}\, \pi^6 \,.
$$
This is one of the numbers computed by A.~Zorich. 

\subsection{Connection with random partitions} 

The quantities \eqref{lala} are, by their
construction, certain sums over all partitions
$\la$. The variable $\hb$ enters this sums
as a weight $e^{-\hb|\la|}$ given to a 
partition $\la$. The leading term of the
$\hb\to+0$ asymptotics, like in \eqref{vmu},
describes certain statistical properties
of random partitions of a very large size.

Some of our formulas admit 
a nice probabilistic interpretation, see
the Appendix. In particular, one can easily
see in our formulas the existence of Vershik's
limit shape of a large random partition, see 
Section \ref{A1} and
also the Gaussian correction to this limit 
shape, see Section \ref{CLT}.  

The point of view of random partitions
also provides a very simple explanation
why something like \eqref{Cq} can never
be modular, see Section \ref{fluc}, which
makes the quasimodularity of \eqref{Cq}
look even more like a miracle.

\subsection{Some open problems}
The space $\cH_1(\mu)$ is sometimes disconnected. The
connected components of this space were described in
\cite{Kontsevich:Zorich:components}. In particular,
there are always at most $3$ components and $\cH_1(\mu)$
is connected when at least one of the $\mu_i$'s is odd.
The knowledge if the volumes of connected components
is important for applications to ergodic theory. 
For small genus, volumes of connected components were
determined by A.~Zorich. Unfortunately, our formulas do 
not separate the connected components. 

Another problem important for applications is to
compute the volumes of similarly defined moduli spaces of
quadratic differentials.
  
\subsection{Acknowledgements}
We would like to thank M.~Kontsevich, H.~Masur and A.~Zorich for 
useful conversations, in particular related to
Proposition~\ref{prop:count:tori:and:volumes}.

\section{Counting ramified covering of a torus}
\label{s2}

\subsection{Basics} 

Let $T$ be a torus and $Z=\{z_1,\dots,z_s\}$ be a collection
of distinct points in $T$. 
Let $\sigma:\Sigma\to T$ be a ramified
covering of $T$ which is unramified outside of $Z$.
All information about $\sigma$ is encoded in the 
monodromy action of the fundamental group $\pi_1(T\setminus Z,*)$
on the fiber over the basepoint $*\in T$
$$
\pi_1(T\setminus Z,*) \to \Aut(\sigma^{-1}(*))\,.
$$
If $\sigma$ is $d$-fold then any labeling of $\sigma^{-1}(*)$
by $1,\dots,d$ produces an isomorphism 
$$
\Aut(\sigma^{-1}(*)) \cong S(d)\,.
$$
Therefore, $d$-fold ramified coverings are in bijection with the orbits
of the $S(d)$-action by conjugation on the set of all homomorphisms 
from $\pi_1(T\setminus Z)$ to $S(d)$
\begin{equation}\label{orbitspace}
\Big\{\textup{$d$-fold coverings}\Big\} = \Hom(\pi_1(T\setminus Z),S(d))
\Big/ S(d) \,.
\end{equation}

Introduce the following notation. 
For any conjugacy classes $C_1,\dots,C_s\subset S(d)$, denote by
$$
H_d(C_1,\dots,C_s) \subset \Hom(\pi_1(T\setminus Z),S(d))
$$ 
those
homomorphisms that send a small loop around $z_i$ into $C_i$ for $i=1,\dots,s$.
This corresponds to fixing the ramification type (namely $C_i$) over the 
points $z_i\in Z$.

A natural way to count the orbits in \eqref{orbitspace} is to weight
any orbit $\sigma$ in \eqref{orbitspace} by $|\Aut(\sigma)|^{-1}$ where
$\Aut(\sigma)$ is a point stabilizer of  $\sigma$, that is, 
the centralizer of the image of $\pi_1(T\setminus Z)$ inside $S(d)$.
Introduce the following weighted number of the $d$-fold coverings
with prescribed monodromy $C_1,\dots,C_s$
\begin{align*}
\Cov_d(C_1,\dots,C_s)&=\sum_{\sigma\in H_d(C_1,\dots,C_s) / S(d)} 
\frac 1{|\Aut(\sigma)|} \\
&= \big|H_d(C_1,\dots,C_s)\big| \Big/ d! \,\,.
\end{align*}

Since the conjugacy classes of $S(d)$ are naturally embedded into conjugacy
classes of any bigger symmetric group, it makes sense to introduce
the following generating function
$$
\Cov(C_1,\dots,C_s)=\sum_{d=0}^\infty q^d \, \Cov_d(C_1,\dots,C_s) \,.
$$

\begin{remark} 
To avoid possible confusion, we point out that our definition
of $\Aut(\sigma)$ does not allow permutations of the marked points
$z_1,\dots,z_s$. This will be important in the next subsection
where we consider the relation between connected and disconnected
coverings. For example, if a covering is a union of two
otherwise identical coverings which are ramified over two
different points of $T$ then this covering \emph{does not} have an extra
$\Z_2$-symmetry. 
\end{remark}

\subsection{Connected and disconnected coverings}

The generating function $\Cov(C_1,\dots,C_s)$ counts all, possibly
disconnected, coverings with given monodromy $C_1,\dots,C_s$. In
particular, $\Cov()$ counts all unramified coverings. 

Under the
correspondence \eqref{orbitspace}, connected components correspond
to orbits of the $\pi_1$ action on $\{1,\dots,d\}$ and unramified
connected components correspond to those orbits 
on which small loops around
the $z_i$'s act trivially. Let 
$$
H_d'(C_1,\dots,C_s) \subset H_d(C_1,\dots,C_s)
$$
be the subset corresponding to coverings without unramified connected
components.

\begin{definition}\label{d1}
 Let $\Cov'(C_1,\dots,C_s)$ be the 
generating function for the coverings without unramified connected
components. In other words,
$$
\Cov'(C_1,\dots,C_s) = \sum_{d=0}^\infty q^d \, \frac{|H_d'(C_1,\dots,C_s)|}{d!} \,.
$$
\end{definition}

\begin{definition} 
Similarly, let $\Covc(C_1,\dots,C_s)$ be the generating function for 
connected coverings.
\end{definition}

\begin{lemma}\label{l1}
$$
\Cov'(C_1,\dots,C_s)= \Cov(C_1,\dots,C_s) \Big/ \Cov() \,.
$$
\end{lemma}
\begin{proof}
This is equivalent to
$$
|H_d(C_1,\dots,C_s)|=\sum_{k=0}^d \binom dk\, |H_k'(C_1,\dots,C_s)| \, |H_{d-k}()|\,,
$$
which is obvious.
\end{proof} 

To simplify the exposition, we shall from now on focus on the case when
$C_i$ has a single nontrivial cycle of length $m_i\in\{2,3,\dots\}$.
The case of more general monodromies presents no extra difficulties
but it will be not needed for the application we have in mind.
Accordingly, we shall write $\Cov(m)$, where $m=(m_1,\dots,m_s)$, in
place of $\Cov(C_1,\dots,C_s)$ and similarly for $\CC_d(m)$. 

The function $\Covc$ can be expressed in terms of functions $\Cov'$ as
follows. 
Recall that a partition $\al$ of a set $S$ is a presentation of the set $S$
as an unordered disjoint union of nonempty subsets
$$
S= \al_1 \sqcup \al_2 \sqcup  \dots \sqcup \al_\ell \,,
$$
which are called the blocks of $\al$. 
The number $\ell=\ell(\al)$ is the length of the partition $\al$.  
We denote by $\Pi_s$ the set of all partitions of $\{1,\dots,s\}$.
Any covering $\sigma\in H_d(m)$ produces a partition $\al=\al(\sigma)\in\Pi_s$ as
follows. Two numbers $i$ and $j$ belong to the same block of $\al$ if and
only if the corresponding ramifications occur on the same connected
component. 

It is clear that the same argument that establishes
Lemma \ref{l1} shows that 
\begin{equation}\label{cov'}
\Cov'(m)= \sum_{\al\in\Pi_s} \prod_{k=1}^{\ell(\al)} \Cov^\circ \left(m_{\al_k}\right)  \,,
\end{equation}
where $m_{\al_k}=\{m_i\}_{i\in\al_k}$\,.

\begin{remark}\label{Moe}
Recall that the set $\Pi_n$ of partitions
of an $n$-element set is partially ordered: if $\al,\beta\in\Pi_n$
we say that $\al<\beta$ if $\al$ is a refinement of 
$\beta$, that is, if the blocks of $\beta$ consist of whole blocks
of $\al$. The maximal element of this poset is the partition
$\hn$ into one block.  The
M\"obius function of the partially ordered set
$\Pi_n$ is well known to satisfy
$$
\textup{M\"obius}\,(\hn,\al)=(-1)^{\ell(\al)-1} (\ell(\al)-1)!\,,
$$
see, for example, Section 3.10.4 of \cite{St}. 
\end{remark}

Applying the M\"obius inversion to \eqref{cov'} 
results in the following: 

\begin{lemma}\label{l2} We have
$$
\Covc(m)=\sum_{\al\in\Pi_s} (-1)^{\ell(\al)-1} (\ell(\al)-1)!  
\prod_{k=1}^{\ell(\al)} \Cov' \left(m_{\al_k}\right) \,,
$$
where $m_{\al_k}=\{m_i\}_{i\in\al_k}$\,.
\end{lemma}

\subsection{Coverings and sums over partitions} 

\begin{definition}
Let $C$ be a conjugacy class in $S(d)$.
Let $f_C$ be  the following function of a partition $\la$ 
$$
f_{C}(\la) = \#C \, \frac{\chi^\la(C)}{\dim\la}\,,
$$
where $\chi^\la$ is the character of the irreducible representation of
$S(d)$ corresponding to the partition $\la$, $\chi^\la(C)$ is its
value on any element of $C$, and $\dim\la=\chi^\la(1)$ is the dimension
of representation $\la$.
\end{definition}

If $C$ is the class of an $m$-cycle we shall write $f_m$ instead of $f_C$. 
Also, note the difference between partitions and partitions of a set. In the 
above definition we have simply partitions whereas in the previous section we used
partitions of a set. 

For the number of ramified coverings, there exists the following
expression in terms of the function $f_C$ which goes back
essentially to Burnside, see Exercise 7 in \S 238 of \cite{Burn}. 
In exactly this form it is presented,
for example, in \cite{D}.   

\begin{proposition} We have  
\label{dw}
$$
\Cov_d(C_1,\dots,C_s) = \sum_{|\la|=d} \prod_{i=1}^s f_{C_i} (\la) \,,
$$
where the sum is over all partitions $\la$ of the number $d$. 
\end{proposition}

It is known that for any conjugacy class $C$ the function $f_C(\la)$ is a
polynomial function of a partition $\la$ in the following sense. 

Let $\Ls(n)$ be the 
algebra of polynomials in $\la_1, \dots,\la_n$ which are symmetric
in the variables $\la_i-i$. This algebra is filtered by the degree
of a polynomial.
Let the algebra $\Ls$ be the projective limit
of these algebras $\Ls:=\varprojlim \Ls(n)$ as filtered algebras with respect to
homomorphisms that set the last variable to $0$. This is the algebra
of \emph{shifted symmetric functions}, see \cite{KO,OO}. By construction,
any $f\in\Ls$ has a well defined degree and can be evaluated at any
partition $\la$. There is the following result, see \cite{KO} and also \cite{OO}

\begin{proposition}[\cite{KO}] We have $f_C\in\Ls$ and the degree of $f_C$ is
the number of non-fixed points of any permutation from $C$. 
\end{proposition}

Various expressions are known for this polynomial; for example, its
expression in the shifted Schur functions is given by the formula
(15.21) in \cite{OO}. 

It is clear that we have  
$$
\Cov(m)=\sum_\la q^{|\la|} \, \prod_i f_{m_i}(\la) \,,
$$
where the sum is over all partitions $\la$. In particular, the
generating function for the unramified coverings is
$$
\Cov()=\sum_\la q^{|\la|}  = (q)_\infty^{-1} \,,
$$
where $(q)_\infty=\prod_{n\ge 1} (1-q^n)$.  

Introduce the following linear functional on the algebra $\Ls$

\begin{definition}\label{d2} For any $F\in\Ls$, set
$$
\lan F\ran_q = (q)_\infty \sum_\la q^{|\la|} \, F(\la) \,.
$$
In particular, $\lan 1 \ran_q=1$. More generally, for $s=1,2,\dots$
consider the following multilinear functional on $\left(\Ls\right)^{\times s}$
$$
\lan F_1 \big| F_2 \big| \dots \big| F_s \ran_q = 
\sum_{\al\in\Pi_s} (-1)^{\ell(\al)-1} (\ell(\al)-1)!  
\prod_{k=1}^{\ell(\al)} \lan \prod_{i\in\al_k} F_i \ran_q 
$$
\end{definition}

In other words, $\lan f\ran_q$ is the expected value of $f$ if
the probability of a partition $\la$ is proportional to $q^{|\la|}$.
In the physical language,  $\lan f \ran_q$ is the Gibbsian
average of $f$ with respect to the ``energy'' function $\la\mapsto |\la|$.
The functional $\lan F_1 \big| F_2 \big| \dots \big| F_s \ran_q$ in the
physical language would correspond to the ``connected'' part of
$\lan F_1 F_2 \cdots F_s \ran_q$. It is no coincidence that it
counts quite precisely the connected coverings. 

Indeed, the following is an immediate corollary of Lemmas 
\ref{l1} and \ref{l2}.

\begin{proposition} We have
\begin{align}
\Cov'(m) &=\lan f_{m_1} f_{m_2}  \cdots  f_{m_s} \ran_q \,, \\
\Covc(m) &= \lan f_{m_1} \big| f_{m_2} \big| \dots \big| f_{m_s} \ran_q 
\label{e2}\,.
\end{align}
\end{proposition}

\subsection{Formula for $n$-point functions}\label{forn}

Our strategy for evaluation of the quantities \eqref{e2} is the following.
By multilinearity, it suffices to compute $\lan F_1 \big| F_2 \big| \dots \big| F_s \ran_q$,
where the $F_i's$ range over any linear basis of the algebra $\Ls$ and then
expand the functions $f_m$ in this linear basis. 

\begin{remark}
One fact supporting such a roundabout approach, aside of the fact that 
it appears to be very difficult to evaluate \eqref{e2} directly, is the 
following. As our choice of the parameter $q$ for the generating function
suggests, the averages $\lran_q$ have some
modular properties. More concretely, they are  quasi-modular, see \cite{BO}
and below. 
It turns out, however, that \eqref{e2} are linear combinations of 
quasi-modular forms of different weights or, in other words, they are 
inhomogeneous elements of the algebra of quasi-modular forms. The other
basis of $\Ls$, which will be introduced momentarily, does have the 
property that $\lran_q$ takes basis vectors to
homogeneous quasi-modular forms. 
\end{remark}

A very convenient linear
basis of the algebra $\Ls$ is formed by monomials in the following
generators
\begin{equation}\label{pk}
p_k(\la)=\sum_{i=0}^\infty \left[(\la_i-i+1/2)^k-(-i+1/2)^k\right] + 
(1-2^{-k}) \zeta(-k) \,. 
\end{equation}
This peculiar expression is in fact a natural $\zeta$-function 
regularization of the divergent sum $\sum_{i=0}^\infty (\la_i-i+1/2)^k$.
More precisely, since $\la_i=0$ for all but finitely many $i$ the first
sum in \eqref{pk} is finite while the second term in \eqref{pk} is the natural
regularization for $\sum_{i=0}^\infty (-i+1/2)^k$.  

\begin{remark}
It is an experimental fact
that the somewhat annoying $\frac12$'s in the definition of $p_k$ are actually
very useful, see \cite{BO,BOO,KO,O,OO}. In other words, it turns out that the so-called 
\emph{modified Frobenius coordinates}, which are the usual Frobenius coordinates
plus $\frac12$ for the half of a diagonal square, are the most
convenient coordinates
on partitions. For example, these $\frac12$'s make
the $p_k$ behave well under the involution $\omega$ in the algebra $\Ls$,
see Section \ref{inv}. 
\end{remark}

It is also convenient to introduce the following generating function
$$
\e\la(x)= \sum_{i} e^{(\la_i-i+1/2)x}  \,.
$$
This sum converges provided $\Re x>0$ and has a simple pole at $x=0$ with
residue $1$. We have (see the formula (0.18) in \cite{BO}) 
\begin{equation}\label{pfrome}
p_k(\la)= k!\, \big[x^k\big] \, \e\la(x) \,,
\end{equation}
where $[x^k]$ denotes the coefficient of $x^k$ in the Laurent series 
expansion about $x=0$. Therefore, all averages of the form
$\lan \prod p_{k_i} \ran_q$ are encoded in the following
generating function. 

\begin{definition} We call the following generating function
$$
F(x_1,\dots,x_n) = \lan \prod \e\la(x_i) \ran_q
$$
the \emph{$n$-point function}. Similarly, we also consider more
general generating functions
\begin{multline*}
F(x_1,\dots,x_i\,\big|\,x_{i+1},\dots,x_j\,\big|\,x_{j+1},\dots\,\big|\,\dots ,x_n) =\\
\lan \e\la(x_1)\cdots\e\la(x_i)\,\big|\,\e\la(x_{i+1}) \cdots \e\la(x_j)\,\big|\,
\e\la(x_{j+1}) \cdots \,\big|\, \cdots \e\la(x_n)\ran_q
\end{multline*}
which we call the \emph{connected functions}. 
\end{definition}

It is clear that the connected functions are, by Definition \ref{d2}, polynomials
in the $n$-point functions.

The following claim follows immediately from \eqref{pfrome}

\begin{proposition}\label{p2} Let $\mu$ be a  multi-index 
$\mu=(\mu_1,\dots,\mu_n)$. We have
$$
\lan p_\mu \ran_q = \mu ! \,\big[x^\mu\big] \, F(x) \,,
$$
where $x=(x_1,\dots,x_n)$ and, as usual, 
$$
p_\mu=\prod_i p_{\mu_i}\,,\quad 
\mu!=\prod_i \mu_i\,,\quad
x^\mu=\prod_i x^{\mu_i} \,.
$$
Similarly,
$$
\lan p_\mu \,|\, p_\nu \,|\,p_\eta \,|\, \dots \ran_q  = \mu! \,\nu! \, \eta! \cdots \, 
\big[x^\mu y^\nu z^\eta \cdots \big] \, F(x\,|\,y\,|\, z\,|\, \dots) \,.
$$
\end{proposition}

\begin{definition} The quantities 
$$
\lan p_\mu \,|\, p_\nu \,|\,p_\eta \,|\, \dots \ran_q\,,
$$
which appear in the above proposition and which will be of primary
interest to us in this paper, will be called \emph{cumulants}. 
\end{definition}

Proposition \ref{p2} is, of course, only useful if one can compute the $n$-point
functions. The $n$-point functions were computed in \cite{BO} (see also \cite{O})
as certain determinants involving theta functions and their
derivatives. Introduce the following odd genus 1
theta function
$$
\Th(x)=\Th_{\frac12,\frac12}(x;q)= 
\sum_{n\in\Z} (-1)^n q^{\frac{(n+\frac12)^2}{2}} e^{(n+\frac12)x} \,.
$$ 
This is the only odd genus 1 theta functions and its precise normalization is not
really important because the formulas will be homogeneous in $\Th$. 
The formula for the $n$-point functions is the following

\begin{theorem}[\cite{BO}]\label{t1}
We have
\begin{equation}\label{npoint}
F(x_1,\dots,x_n) = 
\sum_{\substack{\textup{all $n!$ permutations} 
\\\textup{of $x_1,\dots,x_n$}}}
\frac
{\det\left[
\dfrac{\Th^{(j-i+1)}(x_1+\dots+x_{n-j})}{(j-i+1)!}
\right]_{i,j=1}^n}
{\Th(x_1) \, \Th(x_1+x_2) \cdots \Th(x_1+\dots+x_n)}\,,
\end{equation}
where in the $n!$ summands the $x_i$'s have to be permuted in all
possible ways, $\Th^{(k)}$ stands for the $k$-th derivative of $\Th$,
and by the usual convention that $1/k!=0$ if $k<0$ we do not have 
negative derivatives.  
\end{theorem}

In principle, one can use this formula to give a formula
for the connected functions but it appears to be difficult to
simplify the answer in any attractive manner. However, in the 
$q\to 1$ limit, which corresponds to the limit of coverings of
very large degree, the situation simplifies and useful formulas
for the connected functions become available.

\section{Coverings of large degree and volumes of moduli
spaces}

\subsection{Coverings with automorphisms}\label{auto}
Suppose that a $d$-fold connected covering
$$
\sigma: \Sigma \to T
$$
has a nontrivial automorphism, that is, suppose that there
exists a permutation $h\in S(d)$ which commutes with
the monodromy subgroup $G\subset S(d)$ of $\sigma$. 

Since the sets of fixed points of $h^k$, $k=1,2,\dots$,
are $G$-stable and $G$ is transitive, they must be
either empty or all of $\{1,\dots,d\}$ for any $k$.
It follows that the cycle type of $h$ is of the form
$$
(\,\underbrace{d_1,\dots,d_1}_{\textup{$d_2$ times}}\,)\,,
$$
for some factorization $d=d_1 d_2$.

Let $\Z_{d_1}$ be the cyclic group generated by $h$. We have
the following factorization $\sigma=\sigma''\circ\sigma'$ 
\begin{equation}\label{factorsig}
\sigma: \Sigma \xrightarrow[\textup{\quad$d_1$-fold\quad}]{\sigma'} \Sigma/\Z_{d_1}  
\xrightarrow[\textup{\quad$d_2$-fold\quad}]{\sigma''} T \,.
\end{equation}
Because the group $G$ commutes with $h$,
the size of $d_1$ is bounded in terms of the ramification
type $\mu$ of $\sigma$. 

On the other hand, the genus
of $ \Sigma/\Z_{d_1}$ is strictly less than the genus
of $\Sigma$ by the Riemann-Hurwitz formula and
the number of ramification points of $\sigma''$ is
at most the number of ramification points of $\sigma$.

We will see in the next Section that the number of
connected genus $g$ coverings  of degree $\le D$
with $\ell$ ramification points grows like $D^{2g+\ell-1}$
as $D\to\infty$. 
Hence the number of coverings admitting a factorization 
of the form \eqref{factorsig} grows slower than
the number of all coverings.

In particular, the proportion of those coverings of degree $\le D$
which have nontrivial automorphisms becomes negligible as $D\to\infty$. 

\subsection{Proof of Proposition \ref{prop:count:tori:and:volumes}}\label{pofp}

Recall that $p_1,\dots,p_\ell$ denote the zeros of $\omega$
and $\mu_i$'s are the corresponding multiplicities. Also 
recall that we choose the  basis 
$$
\{ \gamma_1, \dots, \gamma_n \} \subset H_1(\Sigma,\{p_i\},\zed)
$$ 
so that that $\gamma_i$, $i=1,\dots,2g$, form a standard symplectic basis
of $H_1(\Sigma,\zed)$ and 
$$
\partial\gamma_{2g+i}=[p_{i+1}]-[p_1]\,,\quad i=1,\dots,\ell(\mu)-1\,.
$$
We have the following elementary

\begin{lemma}  {\rm (cf.\ \cite{Veech:surfaces})} 
\label{lemma:criterion:for:cover}
Consider  $\phi=\Phi(\Sigma,\omega) \in \C^{\dim \cH(\mu)}$.
We have $\phi_i\in\Z^2$, $i=1,\dots,2g$,
 if and only if the following holds:
\begin{itemize}
\item[{\rm (a)}] there exists a holomorphic map $\sigma$ from $\Sigma$ to the
  standard torus $T = [0,1]^2$,
\item[{\rm (b)}] $\omega = \sigma^{-1}(dz)$, 
\item[{\rm (c)}] $\{p_i\}$ is the set of critical points of $\sigma$, 
\item[{\rm (d)}] the ramification of $\sigma$ at $p_i$ is of the form  $z\mapsto z^{\mu_i+1}$, 
\item[{\rm (e)}] $\sigma(p_{i+1})-\sigma(p_1)=\phi_{2g+i}\,\, \mod \Z^2$,
\item[{\rm (f)}] the degree of $\sigma$ is 
equal to $\area_w(\Sigma) = \frac{\sqrt{-1}}{2}
\int_\Sigma \omega \wedge \overline{\omega}$. 
\end{itemize}
\end{lemma}

\begin{proof}
The sufficiency of the conditions in the lemma is clear.
To prove necessity, define the map $\sigma$ by
\begin{displaymath}
\sigma(z) = \int_p^z \omega \quad \mod  \zed^2 \,, 
\end{displaymath}
where $p\in \Sigma$ is arbitrary. This map $\sigma$ is well defined
because $\int_\gamma \omega \in \zed^2$ for any closed path 
$\gamma \subset  \Sigma$. The required properties of $\sigma$
follow easily from the definitions.
\end{proof}

We note the map $\sigma$ depends only on $(M,\omega)$ and not on the
choice of homology basis. Now we finish the proof of 
Proposition~\ref{prop:count:tori:and:volumes} as follows. 

Choose a vector $\beta\in\C^{\dim \cH(\mu)}$ such that 
$$
\beta_i\in\Z^2,\quad i=1,\dots,2g\,,
\quad \beta_i \ne \beta_j \mod \Z^2\,, \quad i,j>2g,i\ne j\,.
$$ 
Let a set $E\subset \cH_1(\mu)$ lie in the domain of
a coordinate chart $\Phi$ and denote by $C_{D}$ the cone
$$
C_D=\left\{t\Phi(\Sigma,\omega), (\Sigma,\omega)\in E, t\in [0,\sqrt D]
\right\} \subset \C^{\dim \cH(\mu)} \,.
$$ 
By definition of $\nu$ we have 
$$
D^{-\dim \cH(\mu)}\left| C_D \cap (\Z^{2 \dim \cH(\mu)}+\beta) \right| \to 
\vol(C_1)=\nu(E)\,,\quad D\to\infty\,. 
$$
On the other hand, by Lemma \ref{lemma:criterion:for:cover},
every point of the intersection $ C_D \cap (\Z^{2 \dim \cH(\mu)}+\beta)$
corresponds to a covering $\sigma$ of degree $\le D$ with ramification 
type $\mu$. Thus, $\nu(E)$ is the asymptotics of
the number of those covering which correspond to the subset $E$
of the moduli space. 

Now for the whole moduli space $\cH_1(\mu)$,
it follows from the proof of the finiteness of the volume 
in \cite{Masur:interval,Veech:Gauss}
that for every $\epsilon > 0$ there exists a compact subset
$K_\epsilon \subset \cH_1(\mu)$ such that $\nu(K_\epsilon) \ge 
\nu(\cH_1(\mu)) - \epsilon$ and it is easy to
show that $\cH_1(\mu)$ has a rectifiable
boundary. Hence 
\begin{equation}
\nu(\cH_1(\mu)) = 
\lim_{D\to\infty} D^{-\dim_\C \cH(\mu)} \, \sum_{d=1}^D \CC_d(\mu+\vec1) \,,
\end{equation}
as was to be shown.

\subsection{Large degree coverings and $q\to 1$ asymptotics}\label{s14}

Recall that we introduced in \eqref{defvo} the 
constants $\vo(\mu)$ such that
$$
\sum_{d=0}^D \Covc_d(m) \sim  \vo(m)\, \frac{D^{|m|+1}}
{|m|+1}  \,, \quad D\to\infty\,.
$$
where $|m|=\sum m_i$. 
We now observe that $\vo(m)$ is determined
by the leading order asymptotics of 
$$
\Covc(m)=\lan f_{m_1} \big| f_{m_2} \big| \dots \big| f_{m_s} \ran_q 
$$
as $q\to 1$. Namely,
we have the following proposition which 
follows from the elementary power series identity 
$$
\frac1{1-q}\sum_{d=0}^\infty  q^d\, a_d = \sum_{d=0}^\infty  q^d \sum_{k=0}^d a_k \,.
$$

\begin{proposition}\label{voq}
$$
\lan f_{m_1} \big| f_{m_2} \big| \dots \big| f_{m_s} \ran_q
 = \vo(m)\, \frac{|m|!}{(1-q)^{|m|+1}}+ O\left((1-q)^{-|m|}\right)\,,
\quad q\to 1 \,. 
$$
\end{proposition}

The $q\to 1$ asymptotics of the $n$-point functions
and of the connected functions will be 
considered in the next section.

\section{Asymptotics of connected functions.}
\label{s3}
\subsection{Asymptotics of $n$-point functions.} 

It will be convenient to replace the parameter $q$, $|q|<1$, by
a new parameter $\hb$, $\Re \hb>0$, related to $q$ by 
$$
q=e^{-\hb}  \,.
$$
The $q\to 1$ limit corresponds to the $\hb\to+0$ limit and
$\dfrac1{1-q}\sim\dfrac 1\hb$. 
The following proposition describes the behavior of the
$\Th$-function in this limit:

\begin{proposition} We have
\begin{equation}\label{asth}
\frac{\Th(\hb x,e^{-\hb})}{\Th'(0,e^{-\hb})}
 = \hb\, \frac{\sin(\pi x)}{\pi} \, 
\exp\left(\frac{\hb x^2}2\right) \,
\left(1+O\left(e^{-\frac{4\pi^2}{\hb}}\right)\right) 
\end{equation}
as $\hb\to +0$ uniformly in $x$. This asymptotic relation can be differentiated 
any number of times. 
\end{proposition}

\begin{proof} The Jacobi imaginary transformation (see e.g.\ Section 1.9 in
\cite{Mum}) yields
$$
\Th(\hb x,e^{-\hb}) = \sqrt{\frac{2\pi}\hb} \, \exp\left(\frac{\hb x^2}2\right)
\, \Th\left(-2\pi i x, e^{-\frac{4\pi^2}{\hb }}\right) \,.
$$
We have
$$
\Th\left(-2\pi i x, e^{-\frac{4\pi^2}{\hb }}\right) = \sum_{n\in Z} 
(-1)^n \, \exp\left(-\frac{2\pi^2(n+\frac12)^2}{\hb }\right) \,
e^{-2\pi i(n+\frac12)x} \,.
$$
It is obvious that this series, together with all derivatives, is
dominated in the $\hb \to+0$ limit by only two terms, namely
the terms  with $n=0$ and $n=-1$ 
which combine into a multiple of $\sin(\pi x)$. 
All other terms differ by a factor of at least $O\left(e^{-\frac{4\pi^2}{\hb }}\right)$.  
\end{proof}

\begin{remark} As we will see below, all Laurent coefficients of
all connected functions behave asymptotically like powers of $\hb $ as $\hb \to+0$.
Therefore, error terms of the form $\exp\left(-\textup{const}\big/{\hb }\right)$
are completely negligible. 
\end{remark}

We want to introduce an operation $\cA$ of ``taking the asymptotics''
which replaces all  
$\Th$-functions and their derivatives by their asymptotics as $\hb \to+0$.
Since the $n$-point functions \eqref{npoint} and all connected functions
are homogeneous in $\Th$, we can ignore the constant factor $\hb  \, \Th'(0,e^{-\hb })/\pi$.
Let us, therefore, make the following:  

\begin{definition} Introduce the following substitution
operator $\cA$ 
$$
\cA(g)=g\Big|_{
\displaystyle
\Th(\hb x,e^{-\hb }) \mapsto \sin(\pi x) \, \exp(\hb x^2/2)} \,,
$$
where $g$ is any expression containing $\Th$-functions and their
derivatives. 
\end{definition}

In particular, we have 
\begin{equation}\label{Ader}
\cA\left(\Th^{(k)}\left(\hb x\right)\right) = \frac1{\hb^k}\, \frac{d^k}{dx^k} \, 
\sin(\pi x) \, \exp\left(\frac{\hb x^2}2\right) \,,  
\end{equation}
where $k=0,1,2,\dots$\,. 

\begin{definition} Introduce the following \emph{asymptotic $n$-point function}
$$
A(x_1,\dots,x_n) = \cA\left(F(\hb x_1,\dots,\hb x_n) \right) \,.
$$
In other words, this is the result of substituting \eqref{asth} into
the formula for the $n$-point functions and discarding the error terms. 
Similarly, define the \emph{asymptotic connected functions}
$$
A(x_1,\dots\,|\,\dots\,|\,\dots,x_n)=  \cA\left(F(\hb x_1,\dots\,|\,\dots\,|\,\dots,\hb x_n) \right) \,.
$$ 
\end{definition}

Our next goal is to derive a formula for the asymptotic $n$-point function.
We will see that it is considerably more simple than the $n$-point functions
\eqref{npoint}. 

\begin{definition}\label{dS}
Introduce the following function
$$
\cS(x_1,\dots,x_n)=\frac{\pi (x_1+\cdots+x_n)^{n-1}}{\sin(\pi(x_1+\cdots+x_n))} \,.
$$
More generally, given any  partition $\al\in\Pi_n$ 
$$
\{1,\dots,n\} = \al_1\sqcup \dots \al_{\ell(\al)}
$$
set, by definition,
$$
\cS_\al(x_1,\dots,x_n)=\prod_{k=1}^{\ell(\al)} \cS(x_{\al_k}) \,,
$$
where $x_{\al_k}=\{x_i\}_{i\in\al_k}$. 
\end{definition}

\begin{remark}\label{sine series}
Because, eventually, we will be expanding the functions $\cS$ into
Laurent series we recall the following Taylor series 
$$
\frac{\pi x}{\sin(\pi x)} = \sum_{k=0}^\infty (2-2^{-2k+2}) \, \zeta(2k)\, x^{2k}\,. 
$$
\end{remark}

\begin{theorem}\label{t2} We have
\begin{equation}\label{anpoint}
A(x_1,\dots,x_n)= e^{-\frac \hb 2\left(\sum x_i\right)^2}
\sum_{\al\in\Pi_n} \hb^{-\ell(\al)} \, \cS_\al(x_1,\dots,x_n) \,,
\end{equation}
where the summation is over all partitions $\al$ of the set $\{1,\dots,n\}$
and the functions $S_\al$ were defined in Definition \ref{dS}.
\end{theorem}

\begin{remark}\label{n=1} It is clear that
$$
A(x_1)= \exp\left(-\frac {\hb  x_1^2}2\right)\, \frac{\pi}{\hb \, \sin(\pi x_1)} \,,
$$
and, thus, \eqref{anpoint} is satisfied if $n=1$.
\end{remark}

\begin{remark} Recall that the $n$-point functions are, by their
definition, certain averages over the set of all partitions. As
$q\to 1$, larger and larger partitions play an important role
in these averages, so the $q\to 1$ asymptotics of the $n$-point
functions is, in a sense, the study of a very large random partition,
see the Appendix. In particular, the factorization of the leading
order asymptotics 
$$
A(x_1,\dots,x_n)= \hb^{-n}\, \prod_{i=1}^n \frac \pi{\sin(\pi x_i)} + O(\hb^{-n+1})\,, \quad 
\hb \to +0 \,,
$$
corresponds to the existence of Vershik's
limit shape of a typical large partition.
\end{remark}

The proof of Theorem \ref{t2} will be based on a sequence of lemmas.
First, note that the denominators of all summands in \eqref{npoint} have a 
factor of $\Th(x_1+\dots+x_n)$. It is convenient to set, by definition,
$$
\tF(x_1,\dots,x_n) = \Th(x_1+\dots+x_n)\, F(x_1,\dots,x_n) \,.
$$
Similarly, set 
\begin{multline*}
\tA(x_1,\dots,x_n) = \cA\left(\tF(\hb x_1,\dots,\hb x_n)\right)=\\
\sin\left(\pi\left(\sum x_i\right)\right)
 \, e^{\frac \hb2\left(\sum x_i\right)^2} \, 
A(x_1,\dots,x_n) \,.
\end{multline*}
We have the following

\begin{lemma}\label{l21} The function $\tA(x_1,\dots,x_n)$ is a polynomial 
expression 
in $\hb^{-1}$, the variables $x_i$, and cotangents of the form $\cot\left(\pi\sum_{i\in S} x_i\right)$,
where $S$ is a subset of $\{1,\dots,n\}$. The degree of $\tA$ in $\hb^{-1}$ equals
$n$. 
\end{lemma}

\begin{proof} Observe that all $\Th$-functions appear in 
$\tF(\hb x_1,\dots,\hb x_n)$ in the 
following combinations: either they appear in pairs of the form 
$$
\frac{\Th^{(k)}(\hb y)}{\Th(\hb y)} \,, \quad y=\sum_{i\in S} x_i \,, 
$$
where $S$ is a subset of $\{1,\dots,n\}$, or else they appear as
the nullwerts
$$
\Th^{(k)}(0) \,.
$$
It is clear from \eqref{Ader}
that the asymptotics in the either case is a polynomial in $\hb^{-1}$ of degree
$k$ with coefficients involving $y$ and $\cot(y)$. It remains to observe
that in all monomials which appear in the expansion of the determinant
in \eqref{npoint} the orders of the derivatives sum up to $n$. 
\end{proof}
  
It is clear that $A(x_1,\dots,x_n)$ is meromorphic with
at most first order poles at the divisors
$$
\cD_{S,m}=\left\{\sum_{i\in S} x_i = m\right\} \,, \quad S\subset\{1,\dots,n\}\,, \quad m\in\Z \,,
$$
and no other singularities.

\begin{remark}\label{r2} For any nonsingular
point $x=(x_1,\dots,x_n)$,
the asymptotic $n$-point function $A(x)$ 
describes the polynomial in $\hb$ terms in the asymptotics of $F(\hb x)$
as $\hb\to+0$. More generally, since the asymptotics
\eqref{asth} is uniform in $x$, any nonsingular contour integral of $A$ 
represents the asymptotics of the corresponding integral for $F$.
In particular, the residues of $A$ at the divisors $\cD_{S,m}$ are determined by the
corresponding residues of $F$.  
\end{remark}

\begin{lemma} The function $A(x_1,\dots,x_n)$ is regular at the divisors
$\cD_{S,0}$ provided $|S|>1$. At $\cD_{\{1\},0}=\{x_1=0\}$ we have
$$
A(x_1,\dots,x_n)=\frac 1{\hb x_1} A(x_2,\dots,x_n) + \dots 
$$
where dots stand for regular terms.
\end{lemma}
\begin{proof} Follows, as explained in Remark \ref{r2},
from the corresponding facts for $F$, see Section 9 in \cite{BO} or
Section 3 of \cite{O}.
\end{proof}

\begin{lemma} The function $A(x_1,\dots,x_n)$ satisfies the following
difference equation
\begin{multline}\label{diffeq} 
A(x_1-1,\dots,x_n)=-e^{\hb \left(\sum x_i-\frac12\right)} \times \\
\sum_{S=\{i_1,i_2,\dots\}\subset\{2,\dots,n\}} (-1)^{|S|}
A(x_1+x_{i_1}+x_{i_2}+\cdots,\dots, \widehat{x_{i_1}},\dots, \widehat{x_{i_2}}, \dots ) \,,
\end{multline}
where the sum is over all subsets $S$ of $\{2,\dots,n\}$ and hats mean that the
corresponding terms should be omitted. 
\end{lemma}
\begin{proof} Follows
from the difference equation satisfied by  $F$, see Section 8 in \cite{BO} or
Section 3 of \cite{O}.
\end{proof}

\begin{definition}\label{d212} Given a partition $\al\in\Pi_n$ and a 
subset $S\subset\{1,\dots,n\}$, write $S\subset \al$ if
$S$ is a subset of one of the blocks of $\al$.
\end{definition}

\begin{lemma} The right-hand side of \eqref{anpoint} satisfies the same difference 
equation \eqref{diffeq} as $A$ does.
\end{lemma}
\begin{proof}
Observe that the binomial theorem and the definition of the function $\cS$ imply that
\begin{multline*}
\cS(x_1-1,\dots,x_k)= \\
\sum_{S=\{i_1,i_2,\dots\}\subset\{2,\dots,k\}} (-1)^{|S|+1} 
\cS(x_1+x_{i_1}+x_{i_2}+\cdots, \dots, \widehat{x_{i_1}},\dots, \widehat{x_{i_2}}, \dots )\,,
\end{multline*}
where the sum is over all subsets $S$ of $\{2,\dots,k\}$ and hats mean that the
corresponding terms should be omitted. Interchanging the order of summation
in the partition $\al$ and in the subset $S$ one obtains
\begin{multline}
\sum_{\al\in\Pi_n}\hb^{-\ell(\al)}\,\cS_\al(x_1-1,\dots,x_n)= \\
\sum_{\al\in\Pi_n}\hb^{-\ell(\al)} \sum_{\{1\}\cup S \subset \al} (-1)^{|S|+1} 
\cS_\al(x_1+x_{i_1}+x_{i_2}+\cdots, \dots, \widehat{x_{i_1}},\dots )= \\ 
\sum_{S\subset \{2,\dots,n\}} (-1)^{|S|+1}   \sum_{\al'\in\Pi_{n-|S|}} 
\hb^{-\ell(\al')}\, \cS_{\al'}(x_1+x_{i_1}+x_{i_2}+\cdots, \dots)\,,
\end{multline}
where $\al'$ is a partition of the set with $n-|S|$ elements 
which is obtained from
the partition $\al$ by mapping $\{1\} \cup S$ to a point. Note
that $\{1\}\cup S \subset \al$, which according to
Definition \ref{d212} means that $1$ and $S$ belong to the same
block of $\al$, implies $\ell(\al)=\ell(\al')$.  

Now the obvious identity
\begin{equation}\label{e5}
e^{-\frac \hb 2\left(\sum x_i-1\right)^2}=
e^{\hb \left(\sum x_i-\frac12\right)}\,   e^{-\frac \hb 2\left(\sum x_i\right)^2}
\end{equation}
completes the proof. 
\end{proof}

Now we can complete the proof of Theorem \ref{t2}

\begin{proof}[Proof of Theorem \ref{t2}]
By induction on $n$. The case $n=1$ is clear, see Remark \ref{n=1}.

Suppose $n>2$. Denote by $\Aq(x_1,\dots,x_n)$ the right-hand side of \eqref{anpoint}.
We know that $\Aq$ 
satisfies the same difference equation as $A(x_1,\dots,x_n)$ does.
We claim that it also has the same singularities as $A$ does. 

Indeed, $\Aq$ is regular at the divisors $\cD_{S,0}$, $|S|>0$, because
$\cS(x_1,\dots,x_k)$ is regular at $\{x_1+\dots+x_k=0\}$ provided
$k>0$. It is also clear that on $\{x_1=0\}$ we have 
$$
\Aq(x_1,\dots,x_n)=\frac1{\hb x_1} \Aq(x_2,\dots,x_n)+\dots 
$$
and so, by induction hypothesis, $A$ and $\Aq$ have identical
singularities at all divisors $\cD_{S,0}$. Since they also
satisfy the same difference equation, all of their singularities
are identical. 

It follows that the function
\begin{equation}\label{nev}
\sin\left(\pi\left(\sum x_i\right)\right)
 \, e^{\frac \hb 2\left(\sum x_i\right)^2} 
\Big[A(x)-\Aq(x)\Big]
\end{equation}
is regular everywhere.
By the difference equation, the induction hypothesis, and \eqref{e5}
this function is also periodic in all $x_i$'s with period $1$. From Lemma \ref{l21} 
we conclude that \eqref{nev} grows at most polynomially as $\Im x_i \to \infty$
and, therefore, it is a constant. Since both $A$ and $\Aq$ are regular at 
$\{x_1+\dots+x_n=0\}$, the function \eqref{nev} vanishes there. It follows that
it is identically zero. This completes the proof. 
\end{proof}

We conclude this subsection by the following asymptotic version of 
Proposition \ref{p2}. It is clear from Theorem \ref{t2} that the
asymptotic $n$-point function $A(x_1,\dots,x_n)$ can be expanded
into a Laurent series in $x_1,\dots,x_n$ in the neighborhood of
the origin. Same is true about the asymptotic connected 
functions since they are polynomials in the $n$-point functions.
The Laurent coefficients of these connected functions are
responsible for the $\hb \to+0$ asymptotics of the cumulants:

\begin{proposition}\label{p3} We have 
\begin{multline}
\lan p_\mu \,|\, p_\nu \,|\,p_\eta \,|\, \dots \ran_q  =\\
 \hb^{-|\mu|-|\nu|-|\eta|-\dots}
\, \mu! \,\nu! \, \eta! \cdots \, 
\big[x^\mu y^\nu z^\eta \cdots \big] \, A(x\,|\,y\,|\, z\,|\, \dots) + O(\dots)\,,
\end{multline}
where $A$ is the asymptotic connected function, 
$|\mu|=\sum_i \mu_i$, and $O(\dots)$ stands for an error term of the following
type 
$$
O(\dots)=O\left(\frac{e^{-\textup{const}/\hb}}{\hb^{\textup{const}}}\right) \,.
$$
\end{proposition}
\begin{proof} The Laurent coefficients of $A$ are certain contour
integrals and hence by Remark \ref{r2} they represent the asymptotics
of the corresponding coefficients of $F$.
\end{proof}

\begin{definition} Let $\lal\, \cdot \,\ral$ denote the polynomial
in $\hb^{-1}$ part of the asymptotics of $\lran_q$ as $q=e^{-\hb}\to 1$, that
is, the asymptotics of $\lran_q$ without the exponentially small
terms. 
\end{definition}

For example, Proposition \ref{voq} can be restated as
\begin{equation}\label{cmuh}
\lal f_{m_1} \big| f_{m_2} \big| \dots \big| f_{m_s} \ral
 = \vo(m)\, \frac{|m|!}{\hb^{|m|+1}} + \dots
\end{equation}
where dots stand for terms of smaller degree in $\hb^{-1}$. 

\subsection{Asymptotics of the connected functions}

The following notation will be useful in manipulation the
connected functions. Recall that in Remark \ref{Moe} we introduced
a partial ordering on the set $\Pi_n$ of partitions of an $n$-element
set. 

\begin{definition} Let $Q$ be a sequence of functions $Q(x_1,\dots,x_n)$,
where $n=1,2,\dots$. For any  partition  $\al\in\Pi_n$ set, by definition
$$
Q_\al(x_1,\dots,x_n)=\prod_{\textup{blocks $\al_k$}} Q(x_{\al_k}) \,,
$$
where $x_{\al_k}=\{x_i\}_{i\in\al_k}$. Similarly, for any $\al\in\Pi_n$ introduce
the corresponding connected function
\begin{equation*}
Q\left(\,\big|_\al x\right) = 
Q\left(x_{\al_1}\,\big|\,x_{\al_2}\,\big|\,\dots\right)=
\sum_{\beta\ge\al} (-1)^{\ell(\beta)-1} (\ell(\beta)-1)!\, 
Q_\beta (x) \,.
\end{equation*}
It is clear that definition is consistent with Definitions \ref{d2}, 
\ref{dS}. 
\end{definition}

\begin{definition} Given two partitions $\al$ and $\beta$, denote
by $\al\wedge\beta$ the \emph{meet} of $\al$ and $\beta$, that is,
the minimal partition consisting of whole blocks of both $\al$ and
$\beta$. We say that $\al$ and $\beta$ are \emph{transversal} 
and write $\al\perp\beta$ if
$$
\ell(\al)+\ell(\beta)-\ell(\al\wedge\beta)= n\,.
$$
\end{definition}

\begin{remark}\label{transv}
Transversal pairs of partitions are extremal in the
sense that for any $\al,\beta\in\Pi_n$ we have
$$
\ell(\al)+\ell(\beta)-\ell(\al\wedge\beta) \le n \,.
$$
Indeed, any block $\beta_k$ of $\beta$ can intersect at most
$|\beta_k|$ blocks of $\al$ and therefore
$$
\ell(\al)-\ell(\al\wedge\beta) \le \sum_{k=1}^{\ell(\beta)} (|\beta_k|-1)=
n-\ell(\beta) \,.
$$
In other words, $\al\perp\beta$ if $\beta$ bonds the blocks of $\al$ as
effectively as possible. 
\end{remark}

Our goal in this section is to prove a formula for the leading
order asymptotics of connected functions as $\hb\to 0$. 
In other words, we want to compute the term with the minimal
exponent of $\hb $ in the asymptotic connected functions
$$
A\left(\big|_\rho x\right) =
\sum_{\al\ge\rho} (-1)^{\ell(\al)-1} (\ell(\al)-1)!\, 
A_\al (x) \,,
$$
where $x=(x_1,\dots,x_n)$, $\rho\in\Pi_n$, 
and the $A_\al(x)$'s are products of the asymptotic $n$-point functions.
This leading order asymptotics  is described in the following 

\begin{theorem}\label{t3} As $\hb \to+0$ we have 
$$
A\left(\big|_\rho x\right) = 
\hb^{-n+\ell(\rho)-1} \sum_{\al\perp \rho}
 \cS_\al(x) \, \cT_{\al\wedge \rho}(x)+ O\left(\hb^{-n+\ell(\rho)}\right)\,,
$$
where
$$
\cT_{\beta}(x) = (-1)^{\ell(\beta)-1}
\left(\sum  x_i\right)^{\ell(\beta)-2} 
\prod_{\textup{blocks $\beta_k$}} 
\left(\,\sum_{i\in \beta_k} x_i\right)\,.
$$
\end{theorem}

\begin{remark} Observe that if $\ell(\beta)=1$ then $\cT_{\beta}(x)=1$.
\end{remark}

In preparation for the proof of Theorem \ref{t3} we  
introduce the following function
$$
E(x_1,\dots,x_n)=\exp\left(-\frac \hb 2\left(\sum x_i\right)^2\right) \,.
$$
It is clear that
$$
E(x)=1+O(\hb ) \,, \quad \hb \to 0 \,.
$$
The next proposition describes the $\hb \to 0$ asymptotics of the connected
versions of $E$

\begin{proposition}\label{connE}
Let $\rho\in\Pi_n$ be a partition. As
$\hb\to 0$ we have 
$$
E\left(\big|_\rho x\right) = \hb^{\ell(\rho)-1}\, \cT_\rho(x)  + O(\hb^{\ell(\rho)})\,.
$$
\end{proposition}
\begin{proof} Recall that, by definition, 
$$
E\left(\big|_\rho x\right)=\sum_{\al\ge 
\rho} (-1)^{\ell(\al)-1} (\ell(\al)-1)!\, 
E_\al (x) \,.
$$
We have
\begin{equation}\label{e4}
E_\al(x)=\exp\left(-\frac \hb 2\sum x_i^2\right)\, 
\exp\left(-\hb \sum_{\{i\ne j\}\subset \al} x_i x_j \right)\,,
\end{equation}
where, we recall Definition \ref{d212}, $\{i,j\}\subset\al$ means that 
$\{i,j\}$ is a subset of a block of $\al$. 

The first factor in \eqref{e4}
is a common factor for all $\al$. 
The Taylor series expansion of the second factor in \eqref{e4} can
be interpreted as summation over certain graphs $\Gamma$ with
multiple edges
$$
\exp\left(-\hb \sum_{\{i\ne j\}\subset \al} x_i x_j \right)=
\sum_{\Gamma\subset\al}(-\hb )^{\sum m(e)} \prod_{\textup{edges $e=\{i,j\}$}} 
\frac{(x_{i} x_{j})^{m(e)}}{m(e)!} \,,
$$
where $\Gamma\subset\al$
means that $e=\{i,j\}\subset\al$ for any egde $e$ of $\Gamma$, no edges from a 
vertex to itself are allowed, and $m(e)$ is a nonnegative integer, called 
multiplicity, which is assigned to any edge $e$.  

The M\"obius inversion in the partially ordered set $\Pi_n$, see
Remark \ref{Moe}, implies that
$$
E\left(\big|_\rho x\right)=e^{-\frac \hb 2\sum x_i^2} 
\sum_{\textup{$\rho$-connected $\Gamma$}} 
(-\hb )^{\sum m(e)} \prod_{\textup{edges $e=\{i,j\}$}} 
\frac{(x_{i} x_{j})^{m(e)}}{m(e)!} \,,
$$
where $\rho$-connected means that $\Gamma$ becomes connected
after collapsing all blocks of $\rho$ to points, that is,
after passing to the quotient
$$
\{1,\dots,n\} \to  \{1,\dots,n\}\big/\rho\cong\{1,\dots,\ell(\rho)\}\,.
$$

It is now clear that the minimal possible exponent of $\hb $ is $\ell(\rho)-1$
and it is achieved by those graphs $\Gamma$ which have no multiple edges
and project onto spanning trees of $\{1,\dots,\ell(\rho)\}$. That is, 
$$
E\left(\big|_\rho x\right)=(-\hb )^{\ell(\rho)-1} 
\sum_{\substack{\textup{spanning trees}\\\textup{on $\{1,\dots,\ell(\rho)\}$}}} \,\,\,
\prod_{\textup{edges $e=\{k,l\}$}} 
y_k y_l+ O(\hb^{\ell(\rho)})\,,
$$
where $y_k=\sum_{i\in\rho_k} x_i$. It is known, see Problem 3.3.44 in \cite{GJ},
that this sum over spanning trees equals
$$
\sum_{\textup{spanning trees}}=\left(\sum y_k\right)^{\ell(\rho)-2}\, \prod_k y_k\,,
$$ 
which concludes the proof. 
\end{proof}

\begin{remark}\label{spanfor}
 Call a forest with vertices $\{1,\dots,n\}$ a \emph{$\rho$-spanning 
forest} if it has $\ell(\rho)-1$ edges and connects all blocks of $\rho$. 
Equivalently, a forest is $\rho$-spanning if it projects onto a spanning tree
on the quotient $ \{1,\dots,n\}\big/\rho$.

It is clear from the proof of the above proposition that
$$
\cT_\rho= (-1)^{\ell(\rho)-1} \sum_{\substack{\textup{$\rho$-spanning}\\ 
\textup{forests}}} \,\, \prod_{\textup{edges $e=\{i,j\}$}} 
x_i x_j \,.
$$
\end{remark}

\begin{proof}[Proof of Theorem \ref{t3}] By definition, we have
$$
A\left(\big|_\rho x\right) =
\sum_{\beta\ge\rho} (-1)^{\ell(\beta)-1} (\ell(\beta)-1)!\, 
A_\beta (x) \,.
$$
Substituting Theorem \ref{t2} in this sum yields
$$
A\left(\big|_\rho x\right) =
\sum_{\beta\ge\rho} (-1)^{\ell(\beta)-1} (\ell(\beta)-1)!\, 
E_\beta(x) \, \sum_{\al\le\beta} \hb^{-\ell(\al)} \, \cS_\al(x) \,.
$$
Interchanging the order of summation we obtain
\begin{align*}
A\left(\big|_\rho x\right) &=
\sum_{\al}  \hb^{-\ell(\al)} \, \cS_\al(x) 
\sum_{\beta\ge\al\wedge\rho} (-1)^{\ell(\beta)-1} (\ell(\beta)-1)!\, E_\beta(x)\\
&=
\sum_{\al}  \hb^{-\ell(\al)} \, \cS_\al(x) \, E\left(\big|_{\al\wedge\rho} x\right)
\,.
\end{align*}
{}Using Proposition \ref{connE} we conclude that
$$
A\left(\big|_\rho x\right)=\sum_{\al}  \hb^{-\ell(\al)+\ell(\al\wedge\rho)-1} \, \cS_\al(x) \,
\left(\cT_{\al\wedge\rho} + \dots\right)\,,
$$ 
where dots stand for lower order terms. We know from Remark \ref{transv} that
the exponent $-\ell(\al)+\ell(\al\wedge\rho)-1$ takes its minimal value
$-n+\ell(\rho)-1$ precisely when $\al\perp\rho$. This concludes the
proof. 
\end{proof}

\begin{definition} Introduce the following notation for the 
coefficient of $\hb$ in the leading asymptotics of the connected functions
\begin{equation}\label{Alead}
\Al\left(\big|_\rho x\right) = 
\sum_{\al\perp \rho}
\cS_\al(x)\, \cT_{\al\wedge \rho}(x)  \,.
\end{equation}
\end{definition}

It is clear that Proposition \ref{p3} and Theorem \ref{t3} imply 
the formula for the leading asymptotics of the cumulants as $\hb\to+0$

\begin{definition}\label{defweight}
 We call the number $\wt(\mu)=|\mu|+\ell(\mu)$
the weight of a partition $\mu$.
\end{definition}

\begin{theorem}\label{t3'} Let $\mu,\dots,\eta$ be a collection of $s$ partitions.
Then
\begin{equation}
\lan p_\mu \,|\, \dots \,|\, p_\eta \ran_\hb 
=
\frac{\mu! \cdots \eta!\, 
\big[x^{\mu} \cdots z^\eta \big] \, \Al(x\,|\, \dots\,|\,z)}
{\hb^{\wt(\mu)+\dots+\wt(\eta)-s+1}} + \dots\,,
\end{equation}
where $\wt(\mu)=|\mu|+\ell(\mu)$ and dots stand for terms of smaller
degree in $\hb^{-1}$. 
\end{theorem}

We will address the task of actually picking the Laurent coefficients of
$\Al\left(\big|_\rho x\right)$ below in Sections \ref{as cum} and \ref{222}. 
First, we
take a small detour and consider the properties of the weight
function $\wt(\mu)$ which was introduced in Theorem \ref{t3'}

\section{Weight filtration in $\Ls$}
\label{s4}
\subsection{Weight grading and weight filtration} 

The weight function $\wt(\mu)=|\mu|+\ell(\mu)$ introduced in Definition
\ref{defweight} has the following interpretation. 

It is known, see \cite{BO}, and can be seen from the 
formula \eqref{npoint} for the $n$-point functions, that
for any partition $\mu$ 
$$
\lan p_\mu \ran_q \in QM_{\wt(\mu)} \,, 
$$ 
where $QM_*$ is the graded algebra of the quasi-modular form which is
the polynomial algebra in the Eisenstein series $G_k(q)$, $k=2,4,6$.
Therefore, the \emph{weight grading} of $\Ls$ which is defined by assigning
the generators $\{p_k\}$ the weights   
$$
\wt(p_{k})=k+1 \,, \quad k=1,2,\dots\,, 
$$
is very natural in the sense that the linear map $\lan\,\cdot\,\ran_q:\Ls\to QM_*$
preserves it. 
It is clear that 
$$
\lan p_\mu \,|\, \dots \,|\, p_\eta \ran_q \in QM_{\wt(\mu)+\dots+\wt(\eta)}\,.
$$
Proposition \ref{t3'} says that
$$
\lan p_\mu \,|\, \dots \,|\, p_\eta \ran_\hb = \textup{const} \,
\hb^{-\wt(\mu)-\dots-\wt(\eta)+\textup{\# of partitions}-1}+\dots \,.
$$ 
Since we are interested in the coefficient of the lowest power of $\hb$
which is not identically zero by weight considerations, we introduce
the following 

\begin{definition}
We call the filtration of $\Ls$ associated to the 
weight grading the \emph{weight filtration}. 
\end{definition}
 
It is clear that for 
any $g_1,\dots,g_s\in\Ls$ the constant in the expansion
$$
\lan g_1 \,|\, \dots \,|\, g_s \ran_\hb =\textup{const} \,
\hb^{-\sum\wt(g_i)+s-1}+\dots 
$$ 
depends only on the 
top weight terms of $g_1,\dots,g_s$.  

\subsection{Elementary description of the weight filtration} 

In contrast to the weight grading, the weight filtration is 
very easy to describe in completely elementary terms. 

By construction, the algebra $\Ls$ is a projective limit
of the algebras $\Ls(n)$ of shifted symmetric functions 
in $n$ variables. The algebra $\Ls(n)$ is isomorphic
to the algebra of symmetric polynomials in
$$
\xi_i=\la_i-i+\textup{const}\,, \quad i=1,\dots,n \,,
$$
where any constant will do.

It is easy to see that the induced filtration of $\Ls(n)$
is the same as the one obtained by assigning weight $(k+1)$
to the polynomial 
$$
\bar p_k = \sum \xi_i^k \,, \quad k=1,2,\dots \,.
$$ 

Let $\bar m_\mu\in\Ls(n)$ be the monomial symmetric function 
in the $\xi_i$'s, that is, the sum of all monomials which
can be obtained from $\xi^\mu$ by permuting the $\xi_i$'s. 
Recall that the notation $\mu=1^{\vk_1} 2^{\vk_2} 3^{\vk_3} \dots$
means that $\mu$ has $\vk_k$ parts equal to $k$. The following lemma is 
immediate 

\begin{lemma}\label{l31} For any partition $\mu=1^{\vk_1} 2^{\vk_2} 3^{\vk_3} \dots$ we have
$$
\bar p_\mu =\prod \bar p_{\mu_i} =  \vk ! \, m_\mu+\dots\,,
$$
where dots stand for lower weight terms.
\end{lemma} 

\begin{definition} Define the weight of a monomial $\xi^\mu$ by
$
\wt(\xi^\mu)=\wt(\mu)=|\mu|+\ell(\mu)
$
or, in other words, 
$$
\textup{weight}=\textup{degree}+\#\textup{ of variables}\,. 
$$
\end{definition}

It is clear that the $k$-th  subspace of the weight filtration
is spanned by monomials of weight $\le k$. In other words,
we have the following

\begin{proposition} The weight of any shifted symmetric function $g$ is the 
maximum of the weights of all monomials in $g$. 
\end{proposition}

\subsection{Top weight term of $f_k$} 

The purpose of this subsection is to prove the following formula
for the top weight term of $f_k$ 

\begin{theorem}\label{t4} We have
$$
f_k = k^{-1} \sum_{\wt(\lambda)=k+1} \frac{(-k)^{\ell(\lambda)-1}}{\vk!}\,  p_\lambda+\dots\,,
$$
where the sum is over all partitions $\lambda=1^{\vk_1} 2^{\vk_2} 3^{\vk_3} \dots$ of weight
$k+1$ and dots stand for lower weight terms. 
\end{theorem}

\begin{remark}\label{k-1} In fact, the dots in the above formula stand for terms of weight
at most $k-1$ as will be shown in the next subsection. 
In particular, since there are no partitions of weight $1$ we have
\begin{equation}\label{f2p2}
f_2 = \tfrac12\, {p_2} 
\end{equation}
\end{remark}

\begin{proof}
 We can assume that the number of variables $\la_i$ is
finite and equal to $n\gg 0$ 
and switch to the variables $\xi_i=\la_i+n-i$. It is known, see 
Example I.7.7 in \cite{M} , that
\begin{equation*}
f_k=\frac1k \sum_{i=1}^n {(\xi_i\downarrow k)} \,
\prod_{j\ne i} \left(1-\frac{k}{\xi_i-\xi_j}\right) \,, 
\end{equation*}
where $(\xi_i\downarrow k)=\xi_i(\xi_i-1)\cdots(\xi_i-k+1)$. Expand
all fractions in geometric series 
assuming that
$$
|\xi_1|>|\xi_2|>\dots>|\xi_n|\,.
$$
We have
$$
f_k=
\frac1k \sum_{i=1}^n {(\xi_i\downarrow k)} \,
\prod_{j=1}^{i-1} \left(1+k\sum_{l=0}^\infty \frac{\xi_i^l}{\xi_j^{l+1}}\right)
\prod_{j=i+1}^{n} \left(1-k\sum_{l=0}^\infty \frac{\xi_j^l}{\xi_i^{l+1}}\right)
$$
Now let $\mu$ is a partition of weight $k+1$ and 
let us compute the coefficient of $\xi^\mu$ 
in the above expression. 

Observe that only the first summand produces positive
powers of $\xi_1$ and, moreover, the monomials of maximal weight
come from the expansion of 
\begin{equation}\label{ep1}
\xi_1^k \prod_{j=2}^{n} \left(1-k\sum_{l=0}^\infty \frac{\xi_j^l}{\xi_1^{l+1}}\right)
\end{equation}
Clearly, the coefficient of $\xi^\mu$ in the expansion
of \eqref{ep1} equals 
$(-k)^{\ell(\mu)-1}$. By Lemma \ref{l31} this concludes the proof.
\end{proof}

The statement of Theorem \ref{t4} can be rewritten as follows
\begin{equation*}
f_k =- \frac1{k^2} \sum_{\sum_i (i+1) \vk_i = k+1}
 \prod \frac{(-k\, p_{i})^{\vk_i}}{\vk_i!}+\dots \,,
\end{equation*}
where the summation is over all $(\vk_1,\vk_2,\dots)\in\Z_{\ge 0}^\infty$
satisfying the condition
$\sum_i (i+1) \vk_i = k+1$. This can be restated as follows. 

\begin{proposition} We have
$$
f_k = - k^{-2} \, \left[z^{k+1} \right] \, P(z)^k +\dots\,,
$$
where $ \left[z^{k+1} \right]$ stands for the coefficient of $z^{k+1}$,
the dots stand for the lower order terms, and $P(z)$ is the following
generating function
$$
P(z)=\exp\left(-\sum_{i\ge 1} z^{i+1} \, p_i\right) \,.
$$
\end{proposition}

\subsection{Involution and parity in $\Ls$}\label{inv} 

The algebra $\Ls$ has a natural involutive automorphism $\omega$ which
acts as follows
$$
[\omega\cdot f](\la)=f(\la')\,,
$$
where $\la$ is a partition and $\la'$ the dual partition (that is, the 
result of flipping the diagram of $\la$ along the diagonal), see
Section 4 in \cite{OO}. 

For any permutation $g$, we have 
$$
\chi^{\la'}(g)=\sgn(g) \, \chi^{\la}(g) 
$$
and, therefore, 
$$
\omega\cdot f_k = (-1)^{k+1} f_k \,. 
$$
Similarly it can be shown (for example, by expanding 
the statement of Lemma 5.1 in \cite{BO} into a series) that
$$
\omega\cdot p_k = (-1)^{k+1} p_k \,.
$$
It follows that the expansion of $f_k$ in $p_\mu$ contains only
terms of weight 
$$
\wt(\mu)\equiv k+1 \mod 2 \,,
$$
which justifies Remark \ref{k-1}.

\begin{remark}
Note that since $|\la|=|\la'|$ we have
$$
\lan f \ran_q = \lan \omega\cdot f\ran_q
$$
for any $f\in\Ls$. In particular,
$$
\lan p_\mu \ran_q = \lan f_\mu \ran_q= 0 \,, \quad \wt(\mu)\equiv 1 \mod 2\,,
$$
which, of course, makes sense since there are no
quasimodular forms of odd weight. In terms of coverings, this
parity condition just means that the product of monodromies
of all ramifications has to be an even permutation. 
\end{remark}

\section{Asymptotics of cumulants}\label{as cum}
\label{s5}
\subsection{Analog of Wick's formula for cumulants}

Given a multi-index  $m=(m_1,\dots,m_n)$ and a partition
$\rho\in\Pi_n$, we write 
$$
\lal \big|_\rho\, p_m \ral = \lal \prod_{i\in \rho_1} p_{m_i} 
\, \left|  \, \dots \,\left|\, \prod_{i\in \rho_{\ell(\rho)}} p_{m_i} \right. \right. 
\ral \,.
$$
Recall that the Wick formula is a rule to compute expectations of any polynomial
in Gaussian normal variables $\eta_i$ given means $\lan \eta_i \ran$ and covariances 
$\lan \eta_i \,|\, \eta_j \ran$ of these variables. 
Our purpose in this section is to prove a similar rule which reduces the 
computation of any cumulants $\lal \big|_\rho\, p_m \ral$
to computations of the following elementary ones: 

\begin{definition} We call the coefficients 
$\llan m \rran = \llan m_1,\dots,m_n  \rran$ in the expansion
$$
\lal p_{m_1} \,|\, \dots \,|\, p_{m_n}\ral  =
\frac{\llan m\rran}{\hb^{|m|+1}} +
\dots \,,
$$
the \emph{elementary cumulants}. 
\end{definition}

To state the analog of the Wick rule we need the following: 

\begin{definition} Given two partitions $\al,\beta\in\Pi_n$ we say that they are
\emph{complementary} and write $\al\tp\beta$ if $\al\perp\beta$ and 
$\alpha\wedge\beta=\hn$, where $\hn\in\Pi_n$ is the partition into one 
block. In other words, $\al\tp\beta$ if $\beta$ bonds all parts of $\al$ and 
does so using the minimal number of bonds.
\end{definition}

Now we have the following Wick-type formula: 

\begin{theorem}\label{t5} We have 
$$
\lal \big|_\rho\, p_m \ral = \hb^{-\wt(m)+\ell(\rho)-1}\, 
\sum_{\al\tp\rho}\,\,  \prod_{k=1}^{\ell(\al)} 
\llan m_{\al_k}\rran +\dots \,,
$$
where $m_{\al_k} =\{m_i\}_{i\in\al_k}$, $\wt(m)=\sum (m_i+1)$,
and dots stand for lower order terms. 
\end{theorem}

\begin{example} We have
\begin{multline*}
\hb^{a+b+c+d+2}\, \lal p_{a} \,|\, p_{b} \,|\, p_c\, p_d \ral = 
\llan a,b,c\rran\,\llan d\rran+\llan a,b,d\rran\,\llan c\rran\\
+\llan a,c\rran\,\llan b,d\rran+\llan a,d\rran\,\llan b,c\rran +\dots
\,.
\end{multline*}
\end{example}

\begin{example}\label{lead mult}
Note, in particular, that if $\rho=\hn$ then the only
partition complementary to $\rho$ is the partition into 1-element blocks. It
follows that
$$
\lal p_m \ral  = \hb^{-\wt(m)}\, \prod_i \llan m_i\rran +\dots  \,.
$$
This factorization of the leading order asymptotics corresponds to
the limit shape for uniform measure on partitions, see Section \ref{A1}.

Similarly, we have
$$
\lal p_\mu \big| p_\nu \ral =
 \hb^{-\wt(\mu)-\wt(\nu)+1} 
\sum_{k,l} \llan \mu_k,\nu_l \rran \prod_{i\ne k}  \llan \mu_i \rran
\prod_{j\ne l}  \llan \nu_j \rran + \dots \,,
$$
which is the usual Wick's rule for the Gaussian correction to
the limit shape, see Section \ref{CLT}. The covariance matrix
of this Gaussian correction is
$$
\Covar(p_k,p_l)=\hb^{-k-l-1}\, \llan k,l \rran \,.
$$ 
\end{example}

\begin{proof} 
By definition of the cumulants
and of the elementary cumulants, we have 
$$
\lal  p_m \ral = \hb^{-|m|-\ell(\al)}\, 
\sum_{\al}\,\,  \prod_{k=1}^{\ell(\al)} 
\llan m_{\al_k}\rran +\dots \,,
$$
and therefore
$$
\lal \big|_\rho\, p_m \ral = \hb^{-|m|-\ell(\al)}\, 
\sum_{\al\wedge \rho=\hn}\,\,  \prod_{k=1}^{\ell(\al)} 
\llan m_{\al_k}\rran +\dots \,.
$$
By Remark \ref{transv}, for any $\al$ such that $\al\wedge \rho=\hn$,  we have 
$$
|m| + \ell(\al) \le |m| - \ell(\rho) + n + 1 = \wt(m)  - \ell(\rho) + 1\,,
$$
with the equality if and only if $\al\tp\rho$. 
\end{proof}

Theorem \ref{t5} reduces the computation of the asymptotics of cumulants to
the asymptotics $\llan m\rran$ of the elementary cumulants. These numbers will be considered in 
the following subsection.

\subsection{Asymptotics of elementary cumulants}

\begin{definition} We set
$$
\fz(k)=
\begin{cases}
(2-2^{2-k}) \, \zeta(k) &  \textup{$k$ even} \,, \\
0 &  \textup{$k$ odd} \,.
\end{cases}
$$
\end{definition}
By Remark \ref{sine series} this means that
\begin{equation}\label{expsin}
\frac{\pi x}{\sin(\pi x)} = \sum_{k=0}^\infty \fz(k)\, x^{2k}\,. 
\end{equation}

If $\rho$ is a partition into 1-element blocks then any $\al$ is transversal to it
and $\rho\wedge\al=\al$. Therefore,  from Theorem \ref{t3'} and 
\eqref{expsin} we get 
\begin{align*}
\llan m\rran&=m!\, \big[x^{m}\big] \, \sum_{\al}
\cS_\al(x)\, \cT_{\al}(x) \\
&=m!\, \big[x^{m}\big]  \, \sum_{\al}  (-1)^{\ell(\al)-1} \left(\sum x_i\right)^{\ell(\al)-2}
\prod_{k=1}^{\ell(\al)} \sum_{j=0}^\infty \fz(j) \left(\sum_{i\in\al_k} x_i\right)^{j+|\al_k|-1} \,,
\end{align*}
where $m=(m_1,\dots,m_n)$ is a multi-index, the sum is over all partitions $\al\in\Pi_n$.
 
Using the expansion
$$ 
\left(\sum x_i\right)^{\ell(\al)-2}=\sum_{d_1,\dots,d_{\ell(\al)}} 
\binom{\ell(\al)-2}{d_1,\dots,d_{\ell(\al)}} \,
\prod_k \left(\sum_{i\in\al_k} x_i\right)^{d_k}
$$
we obtain the following

\begin{theorem}\label{tem} For any  $m=(m_1,\dots,m_n)$ we have 
\begin{multline}\label{e8}
\llan m\rran=\sum_{\al\in\Pi_n} (-1)^{\ell(\al)-1} (\ell(\al)-2)!  \times \\
\sum_{d} 
(d!)^{-1}
\prod_{k=1}^{\ell(\al)} 
\left|m_{\al_k}\right|!\,\, \fz\left(\left|m_{\al_k}\right| -|\al_k|-d_k+1\right) \,,
\end{multline}
where $\left|m_{\al_k}\right|=\sum_{\al_k} m_i$ and  the summation is over all $\ell(\al)$-tuples
$$
d=(d_1,\dots,d_{\ell(\al)})
$$
of nonnegative integers such that $\sum d_k = \ell(\al)-2$ and 
$$
d_k\equiv 1+\left|m_{\al_k}\right| -|\al_k| \mod 2 \,, \quad k=1,\dots,\ell(\al)\,.
$$ 
The $\al=\hn$ term in \eqref{e8} should be understood  as $|m|! \, \fz(|m|-n+2)$. 
\end{theorem}

\begin{remark}\label{size mu} Given a partition $\mu$ with even parts, write
$$
\fz_\mu=\prod_i \fz(\mu_i) \,.
$$
Observe that all $\fz_\mu$ which appear in  \eqref{e8} satisfy
$$
|\mu|=|m|-n+2 \,.
$$
For any $k$, we have $\fz(k)/\pi^{k} \in \ratls$, and hence
$$
\llan m_1,\dots,m_n \rran / \pi^{|m|-n+2} \in \ratls \,.
$$
\end{remark}

\begin{example}\label{covar} 
In particular, we have
\begin{align}
\llan k \rran &= k! \, \fz(k+1) \,, \\
\llan k,l \rran &= (k+l)!\, \fz(k+l) - k!\, l!\, \fz(k)\, \fz(l) \,,
\end{align}
As already mentioned in Example \ref{lead mult}, these formulas describe
the limit shape of a large random partitions and the covariance
matrix for the  central limit correction to it. 
\end{example}

\begin{remark}
In general, all $\fz_\mu$  which appear in \eqref{e8} are distinct. However,
for small $m$ there may be many like terms and collecting them may lead
to substantial simplifications. In the next section, we will consider
the most interesting of such special cases, namely, the case of $m=(2,\dots,2)$ which
corresponds to the case of simple branched coverings. 
\end{remark}

\section{The case of simple branched coverings}\label{222}
\label{s6}
Consider the case when all ramifications are simple, that is, their 
monodromies only transpose a pair of sheets of the covering. Such coverings
are enumerated by $\lan f_2\,|\,\dots\,|\,f_2\ran_q$. By virtue
of \eqref{f2p2}, the asymptotics $\vo(2,\dots,2)$ of the number of simple coverings
is the following 
$$
\vo(\underbrace{2,\dots,2}_{\textup{$n$ times}}) = \frac1{(2n)! \, 2^{n}} 
\llan 2, \dots,2 \rran \,. 
$$
The aim of this section is to prove the following

\begin{theorem}\label{t7} We have
\begin{equation}\label{st7}
\frac{\vo(2,\dots,2)}{n!} =
\sum_{\substack{\textup{even $\mu$}\\ |\mu|=n+2}}
\frac{ (-1)^{\ell(\mu)-1}}{\vk! \,(2n-\ell(\mu)+2)!}\, 
\left(\prod_{i} (2\mu_i-3)!! \right)\, \fz_\mu\,, 
\end{equation}
where $n$ is the number of $2$'s, 
the summation is over all even partitions 
$\mu=2^{\vk_2} 4^{\vk_4}  6^{\vk_6} \dots$ of the
number $n+2$, and $\vk! = \vk_2!\, \vk_4! \, \cdots $. 
\end{theorem}

The following lemma is well known and elementary to prove

\begin{lemma}\label{l7} For any function $h$ and any $L=1,2,\dots$ we have
$$
\frac1{L!} \left(\sum_{k=1}^\infty h(k) \, \frac{t^k}{k!} \right)^L =
\sum_{n=l}^\infty \frac{t^n}{n!}\, \sum_{\al\in\Pi_n,\,\ell(\al)=L}\, \prod_1^L h\left(|\al_k|\right)\,,
$$
where the summation is over all partitions $\al\in\Pi_n$ which have exactly $L$ parts. 
\end{lemma}
\begin{proof} Follows by extracting terms of degree $L$ in $h$ from the formula
$$
\exp\left(\sum_{k=1}^\infty h(k) \, \frac{t^k}{k!} \right) =
\prod_{k=1}^\infty \exp\left( h(k) \, \frac{t^k}{k!} \right) =
\sum_{n=l}^\infty \frac{t^n}{n!}\, \sum_{\al\in\Pi_n}\, \prod h\left(|\al_k|\right)\,.
$$
\end{proof}

\begin{proof}[Proof of Theorem \ref{t7}]  
The formula \eqref{e8} specializes to 
\begin{multline}\label{e9}
\llan 2,\dots,2 \rran=\sum_{\al\in\Pi_n} (-1)^{\ell(\al)-1} (\ell(\al)-2)!  \times \\
\sum_{d} 
\prod_{k=1}^{\ell(\al)} 
\frac{(2 \left|{\al_k}\right|)!}{d_k!}\,\, \fz\left(|\al_k|-d_k+1\right) \,, 
\end{multline}
where the summation is over $d=(d_1,\dots,d_{\ell(\al)})$ satisfying the
conditions described above. In particular,  $\sum d_i = \ell(\al)-2$. 

Recall that the $\al=\hn$ term in \eqref{e9} is to be understood as
$(2n)!\, \fz(n+2)$. This is in agreement with the coefficient of
$\fz(n+2)$ in \eqref{st7}. Therefore,  in what follows we can 
assume that $\ell(\al)\ge 2$. 

We know from Remark \ref{size mu} that all $\fz_\mu$ appearing in \eqref{e9} satisfy
$$
|\mu|=n+2 \,.
$$
Let $\mu=2^{\vk_2} 4^{\vk_4} 6^{\vk_6}$ be one such partition
and pick the coefficient of $\fz_\mu$ in 
\eqref{e9}. This means that of $\ell(\al)$ blocks of $\al$
we have to chose $\vk_2$ blocks for which
we take $d_k=|\al_k|-1$  so that to produce the factor of $\fz(2)^{\vk_2}$. After that, we
select $\vk_4$ blocks of $\al$ for which we take $d_k=|\al_k|-3$, and so on. In the
remaining 
$$
\vk_0=\ell(\al)-\ell(\mu)
$$ 
parts of $\al$ we take $d_k=|\al_k|+1$ which
results in the factor $\fz(0)=1$.  This can be imagined as painting the parts of $\al$
into different colors which we call ``0'', ``2'', ``4'' etc. 

Observe that the summands in \eqref{e9} depend not on the actual 
partition $\al$ but rather on the sizes of blocks of a given color. For any color
$s=0,2,4,\dots$ we can use Lemma \ref{l7} with $h_s(k)=\frac{(2k)!}{(k-s+1)!}$ 
and this yields the following formula
\begin{multline}\label{e10}
\big[\fz_\mu\big] \llan 2,\dots,2\rran={n!}\,\big[t^n\big] \,\sum_{l=2}^\infty (-1)^{l-1} (l-2)!  \times\\
\prod_{s=0,2,\dots} \frac1{\vk_s!} 
\left(\sum_{k=1}^\infty \frac{(2k)!}{(k-s+1)!\, k!} \, t^k \right)^{\vk_s} \,,
\end{multline}
where $\vk_0=l-\ell(\mu)$ is the only one of the $\vk_i$'s that depends on $l$. 

We have
$$
\sum_{k=0}^\infty \binom{2k}{k} \, t^k  = \frac{1}{\sqrt{1-4t}}
$$
and, therefore, for $s=2,4,6,\dots$ we obtain
\begin{multline}
\sum_{k=1}^\infty \frac{(2k)!}{(k-s+1)!\, k!} \, t^k = t^{s-1} \, \frac{d^{s-1}}{dt^{s-1}}
\frac{1}{\sqrt{1-4t}} \\
 = 2^{s-1}\, (2s-3)!! \,  \frac{t^{s-1}}{(1-4t)^{s-1/2}}
\end{multline}
For $s=0$, introduce the following notation  
$$
H=\sum_{k=1}^\infty \frac{(2k)!}{(k+1)!\, k!} \, t^k = \frac{1-\sqrt{1-4t}}{2t}-1 \,.
$$
We compute
\begin{align*}
\sum_{l=2}^\infty (-1)^{l-1} \frac{(l-2)!}{(l-\ell(\mu))!}\,  
H^{l-\ell(\mu)} &= (-1)^{\ell(\mu)-1} \, (\ell(\mu)-2) !\, (1+H)^{1-\ell(\mu)}\\
&= (-1)^{\ell(\mu)-1} \, (\ell(\mu)-2) ! \frac{2^{\ell(\mu)-1}\, t^{\ell(\mu)-1}}
{\left(1-\sqrt{1-4t}\right)^{\ell(\mu)-1}}\,.
\end{align*}
Putting it all together using the equalities
$$
\sum s\, \vk_s = |\mu|=n+2\,, \quad \sum \vk_s =\ell(\mu) 
$$
we obtain
\begin{multline*}
\big[\fz_\mu\big] \llan 2,\dots,2\rran=(-1)^{\ell(\mu)-1}\, {n!}\,
\frac{2^{n+1} \, (\ell(\mu)-2)!}{\prod_{s\ge 2} \vk_s !}\, 
\left(\prod_{i} (2\mu_i-3)!! \right)  \times \\
\big[t^n\big] \, \frac{t^{n+1}}{(1-4t)^{n+2-\ell(\mu)/2}\, 
\left(1-\sqrt{1-4t}\right)^{\ell(\mu)-1}}\,,
\end{multline*}
It remains to show that 
$$
\big[t^{-1}] \, \frac{1}{(1-4t)^{n+2-l/2}\, 
\left(1-\sqrt{1-4t}\right)^{l-1}} = \frac12 \, \binom{2n}{l-2}\,.
$$
Recall that the residue of a differential form is independent on the
choice of coordinates.
Using the change of variables 
$z=1-\sqrt{1-4t}$, which implies that $dt=\frac12(1-z)\, dz$, 
we compute
\begin{multline*}
\big[t^{-1}] \, \frac{1}{(1-4t)^{n+2-l/2}\, 
\left(1-\sqrt{1-4t}\right)^{l-1}} =\\
\frac12\, \big[z^{-1}] 
\frac{1}{(1-z)^{2n+3-l}\, 
z^{l-1}} = \frac12 \, \binom{2n}{l-2}\,.
\end{multline*}
This concludes proof.
\end{proof}

\begin{appendix}

\section{Large random partitions} 

\subsection{Leading asymptotics and Vershik's limit shape of a
typical random partition}\label{A1}
The computations we do in this paper can be 
interpreted probabilistically as follows. 
We consider the following probability measure $\fP$
on partitions
$$
\fP(\la)=\frac{e^{-\hb|\la|}}Z\,, 
$$
where $Z$ is the partition function
$$
Z=\sum_{\la} q^{|\la|}=\prod_{n\ge 1}(1-q^n)^{-1}\,, \quad q=e^{-\hb}\,. 
$$
It is a Gibbsian measure with the energy function $\la\mapsto
|\la|$ and inverse temperature $\hb$. 

Our algebra $\Ls$ is
naturally an algebra of functions on partitions and what we
can compute is the polynomial terms 
in the asymptotics of the corresponding expectations
$\lran_\hb$ as $\hb\to+0$. The limit $\hb\to+0$
describes the behavior of random partitions of $N$ as $N\to\infty$
and some properties of $\lran_\hb$ have a nice 
interpretation in these terms.

In particular, we have the factorization of the leading 
order asymptotics, see Example \ref{lead mult}, 
$$
\lal p_m \ral  =  \prod_i \hb^{-m_i-1}\, \llan m_i\rran +\dots  \,,
$$
where the numbers $\llan m_i\rran$ are given by
$$
\llan k \rran  =
\begin{cases}
k! (2-2^{-k+1}) \zeta(k+1)\,, & \textup{$k$ odd} \,,\\
0 \,, & \textup{$k$ even} \,,
\end{cases} 
$$
see Example \ref{covar}. 

It is a general principle that if for some probability $\fm$ measure the map
$$
g  \xrightarrow{\quad \textup{Expectation}\quad } 
\lan g \ran =  \int g \, d\fm
$$
is multiplicative,
 then $\fm$ is a  $\delta$-measure. Indeed, the multiplicativity
implies $\textup{Var}(g)=\lan g^2 \ran - \lan g \ran^2=0$ for any $g$. 

Thus, the multiplicativity of the leading order of
the asymptotics reveals the existence of the limit shape of a typical large
partition.

This limit shape is, of course,  the Vershik's limit shape of a typical
large random partition \cite{V} as we shall see momentarily.
Consider the uniform measure on the set
of all partitions of $N$. Then, as $N\to\infty$, the diagrams of typical partitions
are concentrated, see \cite{V} for more precise statements, in
the neighborhood of the shape bounded by the following curve
$$
\exp\left(-\sqrt{\frac{\zeta(2)}{N}} \, x\right) +
\exp\left(-\sqrt{\frac{\zeta(2)}{N}} \, y\right) =1  \,.
$$
The expected size $|\la|$ of a partition $\la$ with respect to
the measure $\fP$ is 
$$
\lal p_1 \ral=\frac{\zeta(2)}{\hb^2}+\dots \,, \quad \hb\to +0\,,
$$
and its variance is
$$
\textup{Var}(|\la|) = \lal p_1 \,|\,p_1 \ral = \frac{\llan 1,1 \rran}{\hb^3}+
\dots =   \frac{\pi^2}{3\hb^3}+\dots 
=o\left(\lal p_1 \ral^2\right)\,.
$$
Therefore, $\fP$-typical partitions have size 
$\approx\zeta(2)/\hb^2$ and hence 
are concentrated in the neighborhood of the shape bounded by the following
curve 
$$
\Upsilon=\Big\{e^{-\hb x} + e^{-\hb y} =1 \Big\} \,.
$$
We will now check that if a partition $\la$ is close to $\Upsilon$ then 
the $p_k(\la)$ is close to $\lal p_k \ral$. Informally, this can
be stated as follows
$$
\lal p_k \ral = p_k\left(\Upsilon\right)+\dots \,.
$$

By performing the summation along the rows of $\la$, one easily checks
that
$$
p_k(\la)=k \sum_{(i,j)\in\la} (j-i)^{k-1} + \dots \,,
$$
where the summation is over all squares $(i,j)$ in the diagram of $\la$
and dots stand for a linear combination of $p_i(\la)$ with $i<k$. 
If $\la$ is close to $\Upsilon$ then
$$
\big|\{(i,j)\in\la, j-i=m\}\big| \approx \hb^{-1} \ln\left(1+e^{-\hb|m|}\right)\,,
$$
and therefore
$$
p_k(\la)\sim \frac{k}{\hb} \sum_{m\in\Z} m^{k-1} \, \ln\left(1+e^{-\hb|m|}\right) \sim
\frac{k}{\hb^{k+1}} \int_{-\infty}^\infty u^{k-1} \ln\left(1+e^{-|u|}\right) \, du \,.
$$
The last integral obviously vanishes if $k$ is even and for $k$ odd it
 is twice the
value of the following Mellin transform
$$
\int_0^\infty u^{s-1} \ln\left(1+e^{-u}\right) \,du =  
\left (1-{2}^{-s}\right )\Gamma (s)\,\zeta (s+1) \,.
$$
Thus, 
we see that indeed $\lal p_k \ral = p_k\left(\Upsilon\right)$+\dots. 

\subsection{Quasimodularity and limit shape fluctuations}\label{fluc}

We recall that averages $\lan p_\mu \ran_q$ are quasimodular forms
in $q$ of weight $\wt(\mu)=|\mu|+\ell(\mu)$. They are only quasi-modular
and \emph{not modular}. This may look like an unfortunate circumstance
from some other points of view, but is actually very natural from the 
point of view of random partitions. For if $\lan p_\mu \ran_q$
 \emph{were} modular
that would mean that the limit shape $\Upsilon$ is incredibly rigid
in the sense that fluctuations of random partitions around $\Upsilon$
would be very, very small. 

For example, if  $\lan p_\mu \ran_q$ were modular then, because there
is only one empty partition, the variance  
$$
\textup{Var}(p_\mu)= \lan p_\mu \,|\, p_\mu \ran_q
$$
would have no $q^0$ term and hence would 
be a \emph{modular cusp form}. Consequently, we would have 
$$
\textup{Var}(p_\mu) = O\left(\frac{e^{-4\pi^2/\hb}}{\hb^{\textup{const}}}\right)
\,, \quad \hb\to +0\,. 
$$
In other words, not only this variance would not grow (compare
this to $\lan p_\mu \ran_q \propto \hb^{-\wt(\mu)}$), but it actually
would decay to $0$ faster than any power of the parameter $\hb$.

In real life, of course, we have
$$
\textup{Var}(p_\mu) \propto \hb^{-2\wt(\mu)+1}\,,
$$
and we have similar power laws for all other cumulants,
see Theorem \ref{t5}. They describe global fluctuation of 
a random partition about the limit shape $\Upsilon$ to
all orders in $\hb$.

It is curious to notice that, on the level of formulas,
these fluctuations are there 
ultimately because the expectations 
$\lan p_\mu \ran_q$ involve the weight
2 Eisenstein series 
$$
G_2(q)=-\frac{1}{24} + \sum_{n=1}^\infty \left(\sum_{d|n} d\right) \, q^n\,,
$$
which has a two-term modular transformation law
$$
G_2\left(e^{-\hb}\right) = -\frac{4\pi^2}{\hb^2} \, G_2\left(e^{-4\pi^2/\hb}\right)-
\frac 1{2\hb}
$$
and, hence, a two-term polynomial asymptotics as $\hb\to+0$
$$
G_2\left(e^{-\hb}\right) = - \frac{\pi^2}{6\hb^2} -\frac 1{2\hb} + O\left(\frac{e^{-4\pi^2/\hb}}{\hb^2}
\right) \,.
$$
 
\subsection{Central limit theorem}\label{CLT} 
Let us center and scale the random variables $p_k$,
that is, introduce the variables
$$
\tip_k = \hb^{k+1/2} \, \left(p_k - \lal p_k \ral\right) \,.
$$ 
We have $\lal \tip_k \ral=0$ and 
\begin{equation}\label{Covpp}
\Covar\left(\tip_k,\tip_l\right) = \llan k,l \rran + O(\hb)\,,
\end{equation}
as $\hb\to+0$, where, see Example \ref{covar}, 
$$
\llan k,l \rran = (k+l)!\, \fz(k+l) - k!\, l!\, \fz(k)\, \fz(l) \,.
$$
Here
$$
\fz(k)=
\begin{cases}
(2-2^{2-k}) \, \zeta(k) &  \textup{$k$ even} \,, \\
0 &  \textup{$k$ odd} \,.
\end{cases}
$$ 
It follows from our Wick formula, see Theorem \ref{t5},
 that the leading asymptotics of averages of the form 
$$
\lal \tip_{m_1} \cdots \tip_{m_n} \ral
$$
is given by the usual Wick rule with the covariance matrix
\eqref{Covpp}. Hence the variables $\tip_k$ are
asymptotically Gaussian normal with mean zero and
covariance \eqref{Covpp}. They describe the Gaussian
fluctuation of a typical partition of $N$ around its
limit shape. 

Similar central limit theorems are known in the literature
for partitions into distinct parts \cite{VFY} and for
the Plancherel measure on partitions \cite{Ke1}. 

Further terms in the asymptotics of $\lran_\hb$ may
not have such a transparent probabilistic interpretation. 

\subsection{Correlation functions and $n$-point functions}

For $y\in\Z+\frac12$,
consider the following function of a partition $\la$
$$
\delta_{y}(\la)
=
\begin{cases} 1\,, & y\in\{\la_i-i+\frac12\}\,,\\
0\,, & \textup{otherwise} \,.
\end{cases}
$$
Then the averages of the form
$$
\lan \delta_{y_1}(\la) \cdots \delta_{y_k}(\la) \ran_q
$$
represent the probability to find the numbers
$y_1,\dots,y_k$ among the numbers $\{\la_i-i+\frac12\}$.
In other words, they are the \emph{correlation functions}
for the 0/1 random process defined on $\Z+\frac12$ by
$y\mapsto \delta_y(\la)$. 

We can write the function $\e\la$ considered in 
Section \ref{forn} in the following form
$$
\e\la(\xi)= \sum_{i} e^{(\la_i-i+1/2)\xi} =
\sum_{y} e^{\xi y}\, \delta_y(\la)  \,.
$$
Therefore, we have, for example,
$$
F(\xi_1|\xi_2) = \lan \e\la(\xi_1) \big | \e\la(\xi_2) \ran_q 
=\sum_{y_1,y_2} e^{\xi_1 y_1 + \xi_2 y_2} \, 
\lan  \delta_{y_1}(\la)  \big | \delta_{y_2}(\la)  \ran_q \,.
$$ 
In other words, the $n$-point functions and
the connected functions are Laplace transforms 
of the correlation functions and their connected
analogs. 

The correlation functions have a nice 
integral representation, see \cite{O}, from
which using the Laplace method one can,
in principle, derive their asymptotics.
This is another possible approach to the
asymptotics of the $n$-point functions.

In particular, in the leading order of
the asymptotics, the correlation functions
factorize, which means that the local shape
of a large random partitions is a random walk.
This factorization is also reflected 
in the leading order factorization of the correlation
functions.

This triviality of the local shape represents
a striking contrast to the situation with the
Plancherel measure  \cite{BOO}, where the local
properties are nontrivial and interesting. Conversely, the
global properties of the Plancherel measure are
quite simple, whereas in our situation their
behavior is rather involved.

\end{appendix}

\end{document}